\documentclass[11pt]{article}

\usepackage{fullpage,latexsym,amsmath,amsfonts,amscd,amsthm,upref,graphics,epsfig}

\newtheorem{numbered}{chapter}[section]
\newtheorem{theorem}[numbered]{Theorem}

\newtheorem{nothing*}[numbered]{}
\newtheorem{remark}[numbered]{Remark}

\newtheorem{definition}[numbered]{Definition}

\newtheorem{proposition}[numbered]{Proposition}

\newtheorem{example}[numbered]{Example}

\numberwithin{equation}{section} \numberwithin{table}{section}
\numberwithin{figure}{section}

\title{Forced Symmetry Breaking from $SO(3)$ to $SO(2)$ for
Rotating Waves on the Sphere}
\author{Adela N. Comanici\\
\emph{\small{Department of Mathematics, University of Houston}}\\
\emph{\small{651 Philip G. Hoffman Hall}}\\
\emph{\small{Houston, Texas, 77204-3008, USA}}\\
\emph{\small{adela@math.uh.edu}}}

\begin{document}
\maketitle

\begin{abstract}
In this article, we consider a small $SO(2)$-equivariant perturbation of a reaction-diffusion system on the sphere,
which is equivariant with respect to the group $SO(3)$ of all rigid rotations. We consider a normally hyperbolic
$SO(3)$-group orbit of a rotating wave on the sphere that persists to a normally hyperbolic $SO(2)$-invariant
manifold $M(\varepsilon)$. We investigate the effects of this forced symmetry breaking by studying the perturbed
dynamics induced on $M(\varepsilon)$ by the above reaction-diffusion system. We prove that depending on the frequency
vectors of the rotating waves that form the relative equilibrium $SO(3)u_{0}$, these rotating waves will give
$SO(2)$-orbits of rotating waves or $SO(2)$-orbits of modulated rotating waves (if some transversality conditions hold).
The orbital stability of these solutions is established as well.
\end{abstract}

\section{Introduction}

In mathematics and physics the phrase "symmetry breaking" has distinct meanings. The first refers to the frequently
observed phenomenon that a configuration of a physical system satisfying a law (a set of equations) which is invariant
under a group of transformations, may itself be invariant under a subgroup of this group. This is referred to as
\emph{spontaneous symmetry breaking}. The second meaning refers to the problem of explicitly adding symmetry breaking
terms to the equations which describe the system. This is called \emph{induced} or \emph{forced symmetry breaking}.
In another words, \emph{forced symmetry breaking} means the systematic study of an equivariant system of differential
equations which is perturbed slightly so that it loses its symmetry properties partially or completely.\\
The main motivation of this article is the presence of spiral waves in excitable media, especially in cardiac tissue.
Spiral waves arise as stable spatio-temporal patterns in various chemical, physical systems and biological systems,
as well as numerical simulations of reaction-diffusion systems on excitable media with various geometries.
Excitable media are extended non-equilibrium systems having a uniform rest state that is linearly stable but susceptible
to finite perturbations. Spiral waves have been observed experimentally, for instance, in catalysis of platinum
surfaces, Belousov-Zhabotinsky chemical reactions \cite{JSW}, Rayleigh-Bernard convection, and the most
important, cardiac tissue \cite{DPSBJ}. Numerical simulations of spiral waves have been done, for example
in \cite{Am, Ba1, GA, YMY, ZM, ZM1}. Winfree \cite{Win} found spiral waves in the Belousov-Zhabotinsky reaction.\\
It is now believed that spiral and scroll waves that appear in the heart muscle can lead to cardiac arrhythmias
(abnormal rhythms in the heart), giving rise to atrial fluttering or ventricular fibrillation. In normal hearts cardiac
arrhythmias are rare, but in diseased hearts cardiac arrhythmias can become more common. For example, if chambers of the
heart become abnormally large, they are susceptible to serious arrhythmias in which waves are believed to circulate in
a fashion that is similar to the circulation of the Belousov-Zhabotinsky waves in a chemical medium. Real human hearts
are enormously complex three-dimensional structures. In this article, we assume that the geometry of the excitable media
is a sphere, but the forced symmetry breaking from $SO(3)$ to $SO(2)$ allows us to consider either the case of the
sphere with a localized inhomogeneity, or the case of a domain with $SO(2)$ symmetry, obtained by slightly deforming
the sphere. In the case of the cardiac tissue, this is clearly an approximation.\\
In the planar case, a rigidly rotating spiral wave is an example of wave pattern rotating around a center and being well
approximated by an Archimedean spiral wave far from rotation center. Near the rotation center, there is a core region of
the spiral wave, where the front of the wave has a tip, whose structure is considered to be the most important in
understanding the behavior of the whole spiral wave. Typical motions which have been observed for planar spiral waves
are rigid rotation, quasi-periodic-meandering and linear meandering. Using an ad hoc model, Barkley \cite{Ba1} was
the first to realize the key importance of the group $SE(2)$ of all planar translations and rotations in describing
the dynamics and bifurcations of planar spiral waves. From a mathematical point of view, for the planar case, rigidly
rotating spiral waves are examples of rotating waves, meandering spiral waves are examples of modulated rotating waves
and linearly drifting spiral waves are examples of modulated travelling waves. The first rigorous mathematical theory
of the planar spiral waves was done by Wulff \cite{Wu1} using Liapunov-Schmidt reduction on scales of Banach spaces.
These results were generalized by Wulff to general non-compact Lie groups in \cite{Wu2}. Comanici \cite{Co1} and
Chan \cite{Cha} have independently expanded this work to the rotation group $G= SO(3)$, using different approaches.
Sandstede, Scheel and Wulff \cite{SSW} proved a finite-dimensional center bundle reduction theorem near a relative
equilibrium $Gu_{0}$ of an infinite-dimensional vector field on a Banach space $X$ on which acts a finite-dimensional
Lie group (not necessarily compact). Using the results of \cite{SSW}, the Hopf bifurcation from one-armed and multi-armed
rotating spiral wave to meandering spiral waves can be studied \cite{GLM1, SSW}. \\
In summary, the Euclidean-equivariant center bundle approach has been remarkably successful at explaining many of
the experimentally observed dynamics and bifurcations of planar spiral waves. However, the Euclidean symmetry is
an approximation, since no physical experiment is infinite in spatial extent.\\
There are indeed experiments on planar spiral waves which exhibit some phenomena which cannot be explained by Euclidean
symmetry alone: boundary drifting \cite{YP, ZM}, spiral anchoring on localized inhomogeneities \cite{DPSBJ, MMV} and
repelling by localized inhomogeneities \cite{MMV}. Boundary drifting has been observed in situations where size of
the spiral wave core is comparable to the size of the spatial domain of the experiment. In this case, the spiral wave
is attracted to the boundary of the domain and then drifts around the boundary in a meandering motion. The spiral waves
with smaller cores do not drift around the boundary: the spiral wave tip goes straight through the boundary without any
impediment \cite{LOPS}. When local inhomogeneities are present, the spiral wave is attracted to these inhomogeneities,
migrate towards them, and then rotate around them. Some spiral waves are also repelled by these inhomogeneities
\cite{MMV}.\\
There are also experiments in which two dynamic attractors, called \emph{entrainment attractor} and \emph{resonance
attractor} are seen, which attractor is observed depends on the initial conditions of the experiment \cite{GZM}.
All of these phenomena break the Euclidean symmetry and as it has been shown in LeBlanc and Wulff \cite{LW}, these
phenomena are generic consequences of imperfect Euclidean symmetry. Specifically, they have studied the effects of
translation symmetry breaking on normally hyperbolic relative equilibria and normally hyperbolic relative periodic orbits
in general systems of Euclidean-equivariant differential equations which undergo a small perturbation that breaks the
translational symmetry, while preserving rotational symmetry.\\
Also, both numerical and experimental work suggest that anisotropy also can lead to certain dynamical states for spiral
waves which are inconsistent with Euclidean symmetry. Anisotropy can lead to phase-locked two-frequency epicycle spiral
motions, and to complicated quasi-periodic meandering patterns which have overall discrete rotational symmetries.
LeBlanc \cite{Le} has investigated the effects the forced rotational symmetry breaking on spiral wave dynamics.\\
All previous results are valid for planar spiral waves. The interest to consider spiral waves on non-planar surfaces is
motivated by the applicability to problems in physiology (cardiology), biology and chemistry. Therefore, the study of
spiral waves by experiments and numerical simulations of reaction-diffusion systems on the sphere and curved surfaces
have recently been undertaken. In the case of spiral waves on a sphere, the dynamics is expected to be quite different
because any spiral wave starting from a rotating center cannot end at a point. The number of tips of a wave front cannot
be odd, and therefore, the dynamics of spiral waves may acquire new features qualitatively different from the planar case.\\
The dynamics of spiral waves in an excitable reaction-diffusion systems on a sphere was numerically investigated by
\cite{Am, GA, ZM, ZM1}, and \cite{YMY} who employ a spectral method using spherical harmonics as basis functions.
Amjadi \cite{Am} treated the case of a periodic oscillating sphere as well. Maselko \cite{Ma}, as well as Maselko and
Showalter, performed experiments with Belousov-Zhabotinsky chemical waves propagating on the surface of a sphere.
The influence of the topological constraints and the inhomogeneity in the excitability on the geometry and dynamics of
spiral waves on a thin spherical shell of excitable media are presented in \cite{DGK}. Also, rigidly rotating waves on
spherical domains have been studied using kinematical theory, but we do not intend to describe this here.\\
In this article we have studied only the trivial isotropy case, because as far as we are aware, there have been observed
no $m$-spiral waves ($m > 1$) on spherical surfaces. Also, in \cite{DGK} it was numerically verified that there is a critical size of
the sphere below which self-sustained spiral waves cannot exist. Therefore, we use a sphere of an arbitrary,
but fixed radius $r$. Guyard and Lauterbach in \cite{GL1, GL2} studied forced symmetry breaking perturbations for
periodic solutions, but their methods applied to the forced symmetry breaking from $SO(3)$ to $SO(2)$ do not give the
results obtained in this article.\\
There are almost no experiments or numerical simulations on spherical domains with localized inhomogeneities and on
approximatively spherical domains, showing phenomena similar to boundary drifting, spiral anchoring at localized
inhomogeneities or repelling from localized in homogeneities. Following the approach for planar spiral waves, this can
be investigated and we intend to present somewhere else, the numerical simulations illustrating the possibilities of
these phenomena for spiral waves on spherical domains with localized inhomogeneities and on approximatively spherical
domains. Since for planar spiral waves, forced symmetry breaking from $SE(2)$ to $SO(2)$ was successful in explaining
phenomena like boundary drifting, spiral anchoring on localized inhomogeneities or repelling by localized inhomogeneities,
we will treat the same problem here, that is the forced symmetry breaking from $SO(3)$ to $SO(2)$ for rotating waves on
the sphere.\\
We consider the group $SO(2)$ as being diffeomorphic with the subgroup of $SO(3)$ defined by
$\{e^{Q\theta}\ |\ \theta \in [0, 2\pi)\}$, where $Q \in so(3)$ such that $\left | Q \right| = 1$. Section \ref{S:Pert_RD}
is concerned with the set-up of the problem of forced symmetry-breaking from $SO(3)$ to $SO(2)$ for a normally hyperbolic
relative equilibrium $SO(3)u_{0}$. This is done by studying $SO(2)$-equivariant reaction-diffusion systems on the sphere,
that are small perturbations of $SO(3)$-equivariant reaction-diffusion systems on the sphere, having a relative
equilibrium $SO(3)u_{0}$ that persists to a normally hyperbolic $SO(2)$-invariant manifold $M(\varepsilon)$, which is
$SO(2)$-equivariant diffeomorphic to $SO(3)$. We recall the fact that the functional-analytical framework can be found
in \cite{Wu2} for the general Lie groups. In Section \ref{S:Perturbation}, we obtain the $SO(2)$-equivariant
finite-dimensional center manifold reduction and the corresponding reduced differential equations for the perturbed
reaction-diffusion system, which include the general form of the $SO(2)$-equivariant perturbation on the center manifold
near $SO(3)u_{0}$, therefore on $M(\varepsilon)$. The orbit space reduction methods are presented in
Section \ref{S:Projection} and we use them to project the $SO(2)$-equivariant perturbed differential equations
on $M(\varepsilon)$ onto the orbit space $SO(3)/SO(2)$, which is diffeomorphic to the unit sphere $\mathbf{S}^{2}$.
Then, we show that the study of the perturbed dynamics (including the orbital stability of the solutions that persist)
on $M(\varepsilon)$ reduces to the analysis of some differential equations on the unit sphere $\mathbf{S}^{2}$.
For $\varepsilon=0$, the dynamics on $M(\varepsilon)$ consist only of rotating waves. Depending on the relation between
the frequency vectors of these rotating waves and $\overrightarrow Q$, we obtain that the rotating waves will project on
the unit sphere $\mathbf{S}^{2}$ either onto two antipodal hyperbolic equilibria or onto periodic solutions.
Using Implicit Function Theorem and the Poincar\'{e} map, we prove the existence of two smooth branches of hyperbolic
equilibria (generically), and the existence of smooth branches of periodic solutions if some transversality conditions
hold, and determine the stability of these solutions. Section \ref{S:Phase_space} presents the effects of forced
symmetry breaking for the rotating waves in the phase space. Namely, depending on the frequency vectors of the rotating
waves that form the relative equilibrium $SO(3)u_{0}$, these rotating waves will give $SO(2)$-orbits of rotating waves or
$SO(2)$-orbits of modulated rotating waves (if some transversality conditions hold). The orbital stability of these
solutions is established as well. All theorems, for which a reference was not given, are proved in
Section \ref{S:Proof_FSB}. Numerical simulations will be presented somewhere else.

\section{Perturbed Reaction-Diffusion Systems on the Sphere}\label{S:Pert_RD}

We consider a perturbed reaction-diffusion system of the form
\begin{equation}\label{E:FSB_RD}
\frac{\partial{u}}{\partial{t}}(t,x)= D\Delta_{S}{u}(t,x)+F(u(t,x),
\lambda)+ \varepsilon G(u(t,x), x, \varepsilon, \lambda) \mbox{ on }
r\mathbf{S^{2}},
\end{equation}
where $r > 0$, $r\mathbf{S^{2}}$ is the sphere or radius $r$, $\mathbf{S^{2}}$ is the unit sphere in $\mathbb{R}^{3}$
and $u= (u_{1}$,$~u_{2}$, \ldots,~$u_{N}) \colon r\mathbf{S^{2}} \rightarrow \mathbb R^{N}$ for $N \geq 1$,
$\Delta_{S}$ is Laplace-Beltrami operator on $r\mathbf{S^{2}}$,
$D= \left(\begin{array}{ccc}
            d_{1}  & \dots  & 0 \\
            \vdots & \ddots & \vdots\\
            0      & \dots  & d_{N}\\
         \end{array}
    \right)$ with $d_{i} > 0$ for $i=1$,~$2$, \ldots,~$N$ are the diffusion coefficients, $F= (F_{1}$,~$F_{2}$, \ldots,
~$F_{N}) \colon \mathbb{R}^{N} \times \mathbb{R} \rightarrow \mathbb{R}^{N}$ and $G= (G_{1}$,~$G_{2}$, \ldots, ~$G_{N})
\colon \mathbb{R}^{N} \times r\mathbf{S^{2}} \times [0,\infty) \times \mathbb{R} \rightarrow \mathbb{R}^{N}$ are
sufficiently smooth functions for $\left | \lambda \right | $, $\varepsilon \geq 0$ small parameters.
In the following sections we suppress the parameter $\lambda$. Therefore, let $F \colon \mathbb{R}^{N} \rightarrow \mathbb{R}^{N}$
and $G \colon \mathbb{R}^{N} \times r\mathbf{S^{2}} \times [0,\infty) \rightarrow \mathbb{R}^{N}$ be sufficiently smooth
functions such that $F(0)=0$ and $G(0, x, \varepsilon)= 0$ for $\varepsilon \geq 0$ small. Using \cite{Wu2},
the reaction-diffusion system (\ref{E:FSB_RD}) defines a sufficiently smooth parameter-dependent local semiflow
$\mathbf{\Phi} $. We recall the notion of linear action of $SO(3)$ on the usual function spaces, that is the function
$T \colon SO(3) \times \mathbf{Y} \rightarrow \mathbf{Y}$, defined by $T(A, u)(x)= A \cdot u(x)= u(A^{-1}x)$,
where $A \in SO(3)$, $u \in \mathbf{Y}$, $x \in r\mathbf{S^{2}}$.\\
We consider the group $SO(2)$ as being diffeomorphic with the subgroup of $SO(3)$ defined by
$\{e^{Q\theta}\ |\ \theta \in [0, 2\pi)\}$, where $Q \in so(3)$ such that $\left | Q \right| = 1$. Thus, we consider that
the Lie algebra $so(2)$ of $SO(2)$ is isomorphic to $\{Qx \mid x \in \mathbb{R} \}$. We will study the reaction-diffusion
system (\ref{E:FSB_RD}) on the function space $\mathbf{Y}=\mathbf{L^{2}}(r\mathbf{S^{2}}, \mathbb{R}^{N})$.\\
Suppose that the superposition operator associated to $G$, denoted by $G \colon \mathbf{Y^{\alpha}} \times [0,\infty)
\rightarrow \mathbf{Y}$, $G(u,\varepsilon)(x)= G(u(x),x,\varepsilon)$ is $SO(2)$-equivariant with respect to the action
$T$ restricted to $SO(2)$, but not $SO(3)$-equivariant with respect to the action $T$ for any
$\varepsilon \geq 0$. Then, the local semiflow $\mathbf{\Phi}$ is $SO(3)$-equivariant for $\varepsilon= 0$, but it is only
$SO(2)$-equivariant for $\varepsilon > 0$ small (see \cite{Co}).
\begin{example}
For any $x \in r\mathbf{S^{2}}$, we denote by $r(x)= d(x, l)^{2}$, where $l$ is the line containing the vector
$\overrightarrow Q$ and $d$ is the Euclidean distance in $\mathbb{R}^{3}$. Let us define $G(u, x)= \widetilde{G}(u, r(x))$,
where $\widetilde{G} \colon \mathbb{R}^{N} \times \mathbb{R} \rightarrow \mathbb{R}^{N}$ is a sufficiently smooth function.
Then the superposition operator associated to $G$ is $SO(2)$-equivariant, but is not $SO(3)$-equivariant in general.
\end{example}
The definitions of relative equilibrium, relative periodic orbits, rotating waves, modulated rotating waves and tip
position function for an $SO(2)$-equivariant or $SO(3)$-equivariant reaction-diffusion system can be found in \cite{Wu2}
for the general Lie groups or in \cite{Co} for the rotation group $G=SO(3)$.\\
Throughout this paper, we use the following property \cite{Ro}: $AXA^{-1}= B$ if and only if
$\overrightarrow B= A \overrightarrow X$, where $A$, $B \in SO(3)$ and $X \in so(3)$. We also use the following notations:
for $X=\left(\begin{array}{ccc}
                  0 & a & -b \\
                  -a & 0 & c \\
                  b & -c & 0\\
             \end{array}
       \right) \in so(3)$, $\overrightarrow X= \left( \begin{array}{c} c\\ b\\ a\\ \end{array} \right)$ and
$\left | X \right | = \left \| \overrightarrow X \right \| = \sqrt{a^{2}+b^{2}+c^{2}}$ and
$\left \| X \right \| = \sqrt{2}\left \| \overrightarrow X \right \|$.

\section{The General Form of the $SO(2)$-Equivariant Perturbation on the Center Manifold}\label{S:Perturbation}

Throughout this paper we suppose that for $\varepsilon= 0$ we have a normally hyperbolic relative equilibrium
$u_{0} \in \mathbf{Y^{\alpha}}$ such that $\Sigma_{u_{0}}= I_{3}$. Let $\mathbf{\Phi}(t, u_{0},0)= e^{X_{0}t}u_{0}$.
We could have presented only Theorem \ref{thm:DE_FSB_RE} here, but we have chosen to discuss the $SO(2)$-equivariant
center manifold reduction in Theorem \ref{thm:CMR_FSB_RE} and Proposition \ref{prop:Pert_RE} for the sake of
completion.
\begin{proposition}[\cite{BLZ, Fi2, HPS}]\label{prop:persistence_RE}
There exists a sufficiently smooth parameter-dependent $SO(2)$-invariant normally hyperbolic manifold $M(\varepsilon)$
for the reaction-diffusion system (\ref{E:FSB_RD}) such that $M(0)= SO(3)u_{0}$, and, for $\varepsilon > 0$ small,
$M(\varepsilon)$ is diffeomorphic to $SO(3)u_{0}$ and this diffeomorphism is $SO(2)$-equivariant.
\end{proposition}
Since $\Sigma_{u_{0}}= I_{3}$, we have that $SO(3)u_{0}$ is diffeomorphic to $SO(3)$, so $M(\varepsilon)$ is
$SO(2)$-equivariant diffeomorphic to $SO(3)$. Using \cite{BLZ, FSSW, HPS, LW, SSW} we obtain:
\begin{theorem}\cite{Co}\label{thm:CMR_FSB_RE}
Let $L$ be the linearization of the right-hand side of (\ref{E:FSB_RD}) with respect to the rotating wave
$\mathbf{\Phi}(t,u_{0},0)= e^{X_{0}t}u_{0}$ in the co-rotating frame, that is \[L= D\Delta_{S}+D_{u}F(u_{0})-X_{0}.\]
Suppose that:
\begin{enumerate}
\item $\sigma(L) \cap \{z \in \mathbb{C} \mid Re \, (z) \geq 0\}$ is a spectral set with spectral projection $P_{*}$,
and dim$(R(P_{*})) < \infty ;
$ \item the semigroup $e^{Lt}$ satisfies $\left | e^{Lt}|_{R(1-P_{*})} \right | \leq Ce^{-\beta _{0} t}$ for some
$\beta _{0} > 0$ and $C > 0$.
\end{enumerate}
Let $V_{*}$ be the orthogonal complement of $T_{u_{0}}(SO(3)u_{0})$ in $E^{cu}= R(P_{*})$. Then, there exists a
sufficiently smooth parameter-dependent center manifold $M^{cu}(\varepsilon)$ of the normally hyperbolic relative
equilibrium $SO(3)u_{0}$ such that $M(\varepsilon) \subset M^{cu}(\varepsilon)$ and $M^{cu}(\varepsilon)$ is
diffeomorphic to $SO(3) \times V_{*}$ for $\varepsilon \geq 0$ small. For $\varepsilon= 0$, $M^{cu}(\varepsilon)$ is
$SO(3)$-invariant, and for $\varepsilon > 0$ small, $M^{cu}(\varepsilon)$ is $SO(2)$-invariant. Furthermore, there exist
sufficiently smooth functions $f_{1} \colon V_{*} \rightarrow so(3)$, $f_{2} \colon V_{*} \rightarrow V_{*}$,
$G^{2}_{pert} \colon SO(3) \times V_{*}\times [0, \infty) \rightarrow V_{*}$ and  a sufficiently smooth vector field
$G^{1}_{pert} \colon SO(3) \times V_{*}\times [0, \infty) \rightarrow T(SO(3))$ with $G^{1}_{pert} (., q, \varepsilon)$
being $SO(2)$-equivariant, $G^{2}_{pert}(., q, \varepsilon)$ being $SO(2)$-invariant such that any solution of
\begin{equation}\label{E:pertdifCM}
\begin{array}{ccc}
\dot{A} &=& Af_{1}(q)+ \varepsilon G^{1}_{pert}(A,q, \varepsilon)\\
\dot{q} &=& f_{2}(q)+ \varepsilon G^{2}_{pert}(A,q, \varepsilon), \\
\end{array}
\end{equation}
on $SO(3) \times V_{*}$ corresponds to a solution of the reaction-diffusion system (\ref{E:FSB_RD}) on
$M^{cu}(\varepsilon)$ under the diffeomorphic identification for $\varepsilon \geq 0$ small. Also, $f_{1}(0)= X_{0}$,
$f_{2}(0)= 0$ and $\sigma(D_{q}f_{2}(0))= \sigma(Q_{*}L|_{V_{*}})$, where $Q_{*}$ is the projection onto $V_{*}$ along
$T_{u_{0}}(SO(3)u_{0})$. Moreover, $0$ is a hyperbolic equilibrium in $V_{*}$ for the second differential equation in
(\ref{E:pertdifCM}) for $\varepsilon=0$.
\end{theorem}
Henceforth, we will identify $M(\varepsilon)$ with $SO(3)$ and $M^{cu}(\varepsilon)$ with $SO(3) \times V_{*}$ for
$\varepsilon \geq 0$ small. Therefore, we will talk about the semiflow $\mathbf{\Phi}$ on $SO(3)$ and on
$SO(3) \times V_{*}$ for $\varepsilon \geq 0$ small.
\begin{proposition}\label{prop:Pert_RE}
The general form of the $SO(2)$-equivariant perturbation on $M^{cu}(\varepsilon)$ given by
Theorem \ref{thm:CMR_FSB_RE} is
\begin{equation}\label{E:pert}
\begin{array}{ccc}
G^{1}_{pert}(A,q, \varepsilon) & = & A k(A,q,\varepsilon),\\
G^{2}_{pert}(A,q, \varepsilon) & = & h(A,q,\varepsilon),\\
\end{array}
\end{equation}
where $A \in SO(3)$, $q \in V_{*}$, $\varepsilon \geq 0$ is small, and $k \colon SO(3)\times V_{*} \times
[0, \varepsilon_{0}) \rightarrow so(3)$ and $h \colon SO(3) \times V_{*} \times [0, \varepsilon_{0}) \rightarrow so(3)$
are sufficiently smooth $SO(2)$-invariant functions for $\varepsilon \in [0,\varepsilon_{0})$ and $q \in V_{*}$ with
$\varepsilon_{0} >0$ small.
\end{proposition}
The proof of Proposition \ref{prop:Pert_RE} is presented in Section \ref{S:Proof_FSB}. Generically, $M(\varepsilon)$
corresponds to a hyperbolic equilibrium $q_{\varepsilon} \in V_{*}$ near $0$ for $\varepsilon \geq 0$ small.
Therefore, we can restrict the study to the $SO(2)$-equivariant dynamics on $SO(3)$, by substituting $q= q_{\varepsilon}$
in the differential equations (\ref{E:pertdifCM}).
\begin{theorem}\label{thm:DE_FSB_RE}
The perturbed differential equations on $M(\varepsilon)$ given by Theorem \ref{thm:CMR_FSB_RE} are
\begin{equation}\label{E:pertdif}
\dot{A}= A \left [X_{0}+\varepsilon g(A, \varepsilon)\right ]
\end{equation}
where $A \in SO(3)$, $\varepsilon \geq 0$ is small and, up to a constant matrix in $so(3)$ that depends sufficiently
smoothly on $\varepsilon$, $g(A, \varepsilon)= k(A, q_{\varepsilon}, \varepsilon)$ for $k$ defined in Proposition
\ref{prop:Pert_RE}. Thus $g(.,\varepsilon)$ is $SO(2)$-invariant for $\varepsilon \geq 0$ small.
\end{theorem}
The proof of Theorem \ref{thm:DE_FSB_RE} is presented in Section \ref{S:Proof_FSB}.

\section{The Analysis of the Projected Dynamics on the Orbit Space $SO(3)/SO(2)$}\label{S:Projection}

\begin{figure}
\begin{center}
\psfig{file=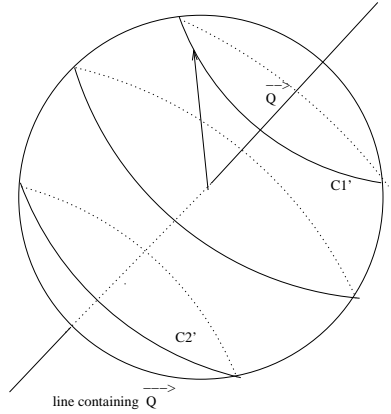,width=2in} \caption{The Sphere $S_{2}$}
\label{figure:figura3}
\end{center}
\end{figure}

\begin{figure}
\begin{center}
\psfig{file=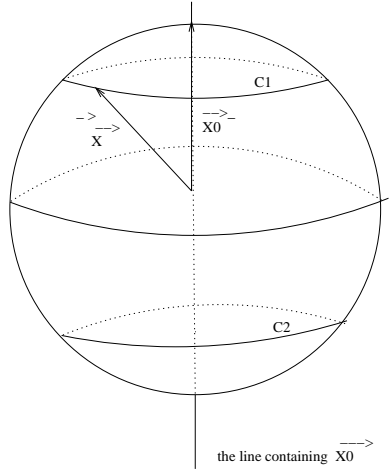,width=2in}
\caption{The Sphere $S_{1}$} \label{figure:figura4}
\end{center}
\end{figure}

Let $x^{0}_{1}$ and $x^{0}_{2}$ be the intersections of the line with direction the vector $\overrightarrow X_{0}$ with
the unit sphere $\mathbf{S}^{2}$; $x^{0}_{1}= \frac{1} {\left | X_{0} \right |}\overrightarrow X_{0}$ and
$x^{0}_{2}= -\frac{1}{\left | X_{0} \right |}\overrightarrow X_{0}$. Let $x^{Q}_{1}$ and  $x^{Q}_{2}$ be the
intersections of the line with direction the vector $\overrightarrow Q$ with the unit sphere  $\mathbf{S}^{2}$;
$x^{Q}_{1}= \overrightarrow Q$ and $x^{Q}_{2}= -\overrightarrow Q$.\\
We consider the left action $\widetilde{\theta} \colon SO(3) \times SO(3) \rightarrow SO(3)$ given by
\begin{equation}\label{E:Action_FSB}
\widetilde{\theta}(A, C)= CA^{-1}.
\end{equation}
We substitute $C= A^{-1}$ in (\ref{E:pertdif}). Since $O_{3}=\dot{I_{3}}=\dot{\widehat{AC}}=\dot{A}C+A\dot{C}$, it
follows that $\dot{C}=-C\dot{A}C$. Therefore, the differential equations $(\ref{E:pertdif})$ are equivalent to the
following differential equations
\begin{equation}\label{E:pertdif_changed}
\dot{C}= -\left [X_{0}+\varepsilon \widetilde{g}(C, \varepsilon)\right ]C,
\end{equation}
where $C \in SO(3)$, $\varepsilon \geq 0$ small and $\widetilde{g}$ is a sufficiently smooth function such that,
for $\varepsilon \geq 0$ small, $\widetilde{g}(C, \varepsilon)=g(C^{-1},\varepsilon)$ is invariant under the restriction
of the action $\widetilde{\theta}$ to the subgroup $SO(2)$ of $SO(3)$, with $g$ defined in Theorem \ref{thm:DE_FSB_RE}.\\
Let us denote the flow associated to the differential equations (\ref{E:pertdif_changed}) by $\mathbf{\Psi}$. It is
clear that $\mathbf{\Psi}(t,C,\varepsilon)=[\mathbf{\Phi}(t,C^{-1},\varepsilon)]^{-1}$. The differential equations
(\ref{E:pertdif_changed}) have the following property:
\begin{itemize}
\item for $\varepsilon= 0$, they are $SO(3)$-equivariant under the action $\widetilde{\theta}$;
\item for $\varepsilon > 0$ small, they are $SO(2)$-equivariant under the restriction of the action $\widetilde{\theta}$
to the subgroup $SO(2)$ of $SO(3)$.
\end{itemize}
We recall the notion of orbit space of $SO(3)$ under the action $\widetilde{\theta}$ (see \cite{Ch}) and some of its
properties that we will use to project the differential equations (\ref{E:pertdif_changed}) from the phase space $SO(3)$
onto the orbit space of $SO(3)$ under the action $\widetilde{\theta}$.
\begin{definition}\label{def:orbit_space}
The quotient of $SO(3)$ under the equivalence relation $A \sim B$ iff $A$ and $B$ are on the same group orbit under
the action $\widetilde{\theta}$ defined in (\ref{E:Action_FSB}) and restricted to $SO(2)$ is called the \emph{orbit space}
of $SO(3)$ under the action $\widetilde{\theta}$ and it is denoted by $SO(3)/SO(2)$. It is the left coset of $SO(3)$
modulo $SO(2)$ under the action (\ref{E:Action_FSB}).
\end{definition}
\begin{theorem}[\cite{Ch, Ko, Wa}]\label{thm:Orbitspace_Manifold}
There is a unique $C^{\infty}$ manifold structure for $SO(3)/SO(2)$ such that the projection $\pi \colon SO(3)
\rightarrow SO(3)/SO(2)$ defined by $\pi(C)= C\cdot SO(2)$ is $C^{\infty}$ and such that there exist local smooth
sections of $SO(3)/SO(2)$ in $SO(3)$; that is, for any $A \cdot SO(2) \in SO(3)/SO(2)$ there exists a local smooth
function $\pi_{s}$ from a neighborhood $U$ of $A \cdot SO(2)$ into $SO(3)$ such that $\pi \circ \pi_{s}= id$ on $U$.
\end{theorem}
\begin{theorem}[\cite{Ch, Ko, Wa}]\label{thm:orbit_space}
The orbit space of $SO(3)$ under the action $\widetilde{\theta}$ is diffeomorphic as a manifold to the unit sphere
$\mathbf{S^{2}}$. The diffeomorphism can be chosen such that
\begin{equation}\label{E:diffeomorphism}
\beta \colon  SO(3)/SO(2) \rightarrow  \mathbf{S^{2}}, \beta(A \cdot SO(2))= Ax_{1}^{Q} \mbox { for any } A \in SO(3).
\end{equation}
\end{theorem}
\begin{definition}\label{def:ProjectedFlow}
The \emph{projected flow} on $SO(3)/SO(2)$ associated to the flow $\mathbf{\Psi}$ is defined by
\begin{equation}\label{E:ProjectedFlow}
\widetilde{\mathbf{\Psi}}(t, \pi(C), \varepsilon)= \pi(\mathbf{\Psi}(t, C, \varepsilon)) \mbox { for any } \varepsilon
\geq 0 \, \mbox{ small }, C \in SO(3) \mbox{ and  } \, t \geq 0.
\end{equation}
\end{definition}
The projected flow on $SO(3)/SO(2)$, $\mathbf{\widetilde{\Psi}}$, is well-defined and sufficiently smooth. From \cite{Ko}
and \cite{Sch}, we know that every sufficiently smooth vector field on $SO(3)/SO(2)$ lifts to a sufficiently
smooth $SO(2)$-equivariant vector field on $SO(3)$, unique up to vector fields tangent to the $SO(2)$-orbits.
\begin{definition}\label{def:DiffProjectedFlow}
Using the diffeomorphism $\beta$ defined in (\ref{thm:orbit_space}), we define the following flow on $\mathbf{S^{2}}$:
\begin{equation}\label{E:DiffProjectedFlow}
\widetilde{\mathbf{\Psi}_{1}}(t, Cx_{1}^{Q}, \varepsilon)= \beta(\widetilde{\mathbf{\Psi}}(t, \pi(C), \varepsilon))
\mbox { for any } \varepsilon \geq 0 \mbox{ small }, C \in SO(3) \mbox{ and } t \geq 0.
\end{equation}
The flows $\mathbf{\widetilde{\Psi}}$ and $\mathbf{\widetilde{\Psi_{1}}}$ are topologically equivalent.
\end{definition}
\begin{theorem}\label{thm:FSB_PROJ_E}
The differential equations given by (\ref{E:pertdif_changed}) projected on the orbit space $SO(3)/SO(2)$ have the
following form (after taking into account Theorem \ref{thm:orbit_space} and Definition \ref{def:DiffProjectedFlow}):
\begin{equation}\label{E:pertdif_projected}
\dot{x}= -\left [X_{0}+\varepsilon g^{S}(x, \varepsilon)\right ]x,
\end{equation}
where $x \in \mathbf{S^{2}}$ and $\varepsilon \geq 0$ is small; that is, the vector field given by
(\ref{E:pertdif_projected}) is the vector field associated to the flow $\widetilde{\mathbf{\Psi}_{1}}(t, Cx_{1}^{Q},
\varepsilon)$. The function $g^{S}(x, \varepsilon)$ is given by
\begin{equation}\label{E:pertonsphere}
g^{S}(x, \varepsilon)= \widetilde{g}(C, \varepsilon), \mbox { where } x= Cx_{1}^{Q} \in \mathbf{S}^{2}, \varepsilon
\geq 0 \mbox{ is small and } \widetilde{g} \mbox{ is defined in } (\ref{E:pertdif_changed}) .
\end{equation}
\end{theorem}
The proof of Theorem \ref{thm:FSB_PROJ_E} is presented in Section \ref{S:Proof_FSB}.
\begin{remark}
\begin{enumerate}
\item For $\varepsilon= 0$, the differential equations (\ref{E:pertdif_projected}) have reflectional symmetry with
respect to the equatorial plane orthogonal to (the line containing) the vector $\overrightarrow X_{0}$.
\item For $\varepsilon > 0$, generically, the differential equations (\ref{E:pertdif_projected}) have no reflectional
symmetry with respect to the equatorial plane orthogonal to (the line containing) the vector $\overrightarrow X_{0}$.
\end{enumerate}
\end{remark}
\begin{proposition}[$\varepsilon= 0$]\label{prop:FSB_Orbit_Phase_E}
\begin{enumerate}
\item For $\varepsilon= 0$ the solutions of the differential equations (\ref{E:pertdif}) are of the form
\begin{equation}
A(t)= A_{0}e^{X_{0}t}, \mbox{ where } A_{0} \in SO(3).
\end{equation}
\item For $\varepsilon= 0$ the solutions of the differential equations (\ref{E:pertdif_projected}) are the equilibrium
points $x^{0}_{1}$ and $x^{0}_{2}$, and the $\frac{2\pi}{\left | X_{0} \right| }$-periodic solutions of the form
$x(t, x_{0})= e^{-X_{0}t}x_{0}$, whenever $x_{0} \neq x^{0}_{1} \mbox{ and } x_{0} \neq x^{0}_{2}$.
\item The correspondence between the solutions of the differential equations (\ref{E:pertdif}) and the
differential equations (\ref{E:pertdif_projected}) is as follows:
\begin{enumerate}
\item The solution $A(t)$ of (\ref{E:pertdif}) corresponding to $A_{0}$ such that $x_{1}^{Q}= A_{0}x^{0}_{1}$ projects
onto the equilibrium $x^{0}_{1}$. \item The solution $A(t)$ of (\ref{E:pertdif}) corresponding to $A_{0}$ such that
$x_{1}^{Q}= A_{0}x^{0}_{2}$ projects onto the equilibrium  $x^{0}_{2}$.
\item The other solutions of (\ref{E:pertdif}) with the initial condition $A_{0}$ projects onto the periodic solutions
of (\ref{E:pertdif_projected}) having the initial condition $x_{0}$ such that $x_{1}^{Q}= A_{0}x_{0}$.
\end{enumerate}
\end{enumerate}
\end{proposition}
The proof of Proposition \ref{prop:FSB_Orbit_Phase_E} is presented in Section \ref{S:Proof_FSB}. For $\varepsilon=0$,
the flow given by the differential equations (\ref{E:pertdif_changed}) is the only $SO(3)$-equivariant flow
that can be obtained when we lift from the differential equations (\ref{E:pertdif_projected}) on $SO(3)/SO(2)$ to $SO(3)$.
\begin{proof}
The differential equations (\ref{E:pertdif_projected}) for $\varepsilon=0$ are $\dot{x}=-X_{0}x$ and its flow is given by
$\widetilde{\mathbf{\Psi}} (t,C \cdot SO(2),0)= e^{-X_{0}t}C \cdot SO(2)$. We get $\pi(\mathbf{\Psi}(t,C,0))=
\widetilde{\mathbf{\Psi}}(t,C \cdot SO(2),0)= e^{-X_{0}t}C \cdot SO(2)$. It follows that $\mathbf{\Psi}(t,C,0)=
e^{-X_{0}t}C e^{\tau(t,C)Q}$. This is $SO(3)$-equivariant under the action $\widetilde{\theta}$ if and only if
$\tau(t,C)=0$ for $t \geq 0$ and $C \in SO(3)$; that is, $\mathbf{\Psi}(t,C,0)$ is the $SO(3)$-equivariant flow given
by differential equations (\ref{E:pertdif_changed}) for $\varepsilon=0$.
\end{proof}
A geometrical interpretation of the dynamics on $M(0)=SO(3)u_{0}$ induced by the reaction-diffusion system
(\ref{E:FSB_RD}) is given below.
\begin{proposition}[$\varepsilon= 0$]\label{thm:RotatingWaves}
We will consider two distinct unit spheres $S_{1}$ and $S_{2}$. Let the points $x_{1}^{Q}$ and $x_{2}^{Q}$ be on
sphere $S_{2}$ as they were defined at the beginning of this section. Also, let the points $x_{1}^{0}$ and $x_{2}^{0}$
be on sphere $S_{1}$ as they were defined at the beginning of this section. For $\varepsilon= 0$, the semiflow
$\mathbf{\Phi}$ on $SO(3)u_{0}$ contains only rotating waves. We have:
\begin{enumerate}
\item The rotating waves for the reaction-diffusion system (\ref{E:FSB_RD}) on $SO(3)u_{0}$ with the unit frequency
vectors $\overrightarrow X$ on a circle around $\overrightarrow Q$ of the unit sphere $S_{2}$ project onto the periodic
solution of (\ref{E:pertdif_projected}), that has as an orbit the circle around $\overrightarrow X_{0}$ of the unit
sphere $S_{1}$ passing through the point $x_{0} \in S_{1}$ such that
$cos(\angle(\vec X,\vec Q))= cos(\angle(\vec X_{0},\vec x_{0}))$.
\item The rotating waves for the reaction-diffusion system (\ref{E:FSB_RD}) on $SO(3)u_{0}$ with the frequency
vector $\left | X_{0} \right | \overrightarrow Q$ project onto $x_{1}^{0}$. The rotating waves for the reaction-diffusion
system (\ref{E:FSB_RD}) on $SO(3)u_{0}$ with the frequency vector $-\left | X_{0} \right | \overrightarrow Q$ project onto
$x_{2}^{0}$. \item Let $C_{1}$ and $C_{2}$ be two distinct circles on the unit sphere $S_{2}$ with the centers on
the line  having the direction vector $\overrightarrow Q$ and such that $d(0, C_{1})= d(0, C_{2})$. Let us denote by
$C_{1}^{'}$ the circle on $S_{1}$ onto which project the rotating waves  on $SO(3)u_{0}$ with the unit frequency vectors
on $C_{1}$. Let us denote by $C_{2}^{'}$ the circle on $S_{1}$ onto which project the rotating waves  on $SO(3)u_{0}$
with the unit frequency vectors on $C_{2}$. Then $d(0, C_{1}^{'})= d(0, C_{2}^{'})$ and $C_{1}^{'} \neq C_{2}^{'}$,
where $d$ is the Euclidean distance in $\mathbb{R}^{3}$.
\end{enumerate}
\end{proposition}
The following two figures illustrate Theorem \ref{thm:RotatingWaves}: The Sphere $S_{2}$ (Figure (\ref{figure:figura3}))
and The Sphere $S_{1}$ (Figure (\ref{figure:figura4})).
The proof of Theorem \ref{thm:RotatingWaves} is presented in Section \ref{S:Proof_FSB}.\\
Let $B \in SO(3)$ such that $Bx_{1}^{0}= (0,0,1)^{T}$. Let $g^{SS}(s, \varepsilon)= Bg^{S}(B^{-1}s, \varepsilon)B^{-1}$,
where $g^{S}$ is defined in $(\ref{E:pertonsphere})$ and for any $s=(x_{1}, x_{2}, x_{3}) \in \mathbf{S}^{2}$ and
$\varepsilon \geq 0$ small,
\[g^{SS}(s,\varepsilon)= g^{SS}(x_{1}, x_{2}, x_{3},\varepsilon)=
                                   \left (
                                       \begin{array}{ccc}
                                      0 & -G_{1}(x_{1}, x_{2}, x_{3},\varepsilon) & G_{2}(x_{1}, x_{2}, x_{3},\varepsilon)\\
                                      G_{1}(x_{1}, x_{2}, x_{3},\varepsilon) & 0  & -G_{3}(x_{1}, x_{2}, x_{3},\varepsilon)\\
                                     -G_{2}(x_{1}, x_{2}, x_{3},\varepsilon) & G_{3}(x_{1}, x_{2}, x_{3},\varepsilon) & 0\\
                                       \end{array}
                                   \right ).
\]
Let \[a(\phi, \theta, \varepsilon)= G_{1}(\sin \phi \cos \theta, \sin \phi \sin \theta, \cos \phi,\varepsilon),\]
\[b(\phi, \theta, \varepsilon)= G_{2}(\sin \phi \cos \theta, \sin \phi \sin \theta, \cos \phi,\varepsilon),\] and
\[c(\phi, \theta, \varepsilon)= G_{3}(\sin \phi \cos \theta, \sin \phi \sin \theta, \cos \phi,\varepsilon),\]
where $\theta \in [0, 2\pi)$, $\phi \in [0, \pi]$, $\varepsilon \geq 0$ small. The analysis of the perturbed differential
equations (\ref{E:pertdif_projected}) on the orbit space $SO(3)/SO(2)$ for $\varepsilon > 0$ small gives us the first
main result of this article:
\begin{theorem}[$\varepsilon > 0$]\label{thm:FSB_DE_Orbit_E}
\begin{enumerate}
\item There exist two sufficiently smooth branches of equilibria for (\ref{E:pertdif_projected}) defined for
$\varepsilon \geq 0$ small, one branch $x^{1}(\varepsilon)$, near $x^{0}_{1}$ and the other one $x^{2}(\varepsilon)$,
near $x^{0}_{2}$, such that $x^{1}(0)= x^{0}_{1}$ and $x^{2}(0)= x^{0}_{2}$.\\
Generically, we have $\frac{\partial{G_{3}}}{\partial{x_{2}}}(0,0,1,0)-\frac{\partial{G_{2}}}{\partial{x_{1}}}(0,0,1,0)
\neq 0$. The stability of the generically hyperbolic equilibria $x^{1}(\varepsilon)$ is as follows:
\begin{enumerate}
\item if $\frac{\partial{G_{3}}}{\partial{x_{2}}}(0,0,1,0)-\frac{\partial{G_{2}}}{\partial{x_{1}}}(0,0,1,0)
< 0$, then $x^{1}(\varepsilon)$ is locally asymptotically stable for the differential equations
(\ref{E:pertdif_projected});
\item if $ \frac{\partial{G_{3}}}{\partial{x_{2}}}(0,0,1,0)-\frac{\partial{G_{2}}}{\partial{x_{1}}}(0,0,1,0)
> 0$, then $x^{1}(\varepsilon)$ is unstable for the differential equations (\ref{E:pertdif_projected}).
\end{enumerate}
Generically, we have $\frac{\partial{G_{3}}}{\partial{x_{2}}}(0,0,-1,0)-\frac{\partial{G_{2}}}{\partial{x_{1}}}(0,0,-1,0)
\neq 0$. The stability of the generically hyperbolic equilibria $x^{2}(\varepsilon)$ is as follows:
\begin{enumerate}
\item if $\frac{\partial{G_{3}}}{\partial{x_{2}}}(0,0,-1,0)-\frac{\partial{G_{2}}}{\partial{x_{1}}}(0,0,-1,0)
< 0$, then $x^{2}(\varepsilon)$ is locally asymptotically stable for the differential equations (\ref{E:pertdif_projected}); \item if $
\frac{\partial{G_{3}}}{\partial{x_{2}}}(0,0,-1,0)-\frac{\partial{G_{2}}}{\partial{x_{1}}}(0,0,-1,0)
> 0$, then $x^{2}(\varepsilon)$ is unstable for the differential equations (\ref{E:pertdif_projected}).
\end{enumerate}
\item Let \[I(\phi)= \int_{0}^{2\pi}\left [b(\phi,t,0)\cos t-c(\phi,t,0)\sin t \right ] \ dt  \] for $\phi \in (0,\pi)$,
and consider\\
$x_{0}= B^{-1}\left (\begin{array}{c}
                             \sin \phi_{0}\cos \theta_{0} \\
                             \sin \phi_{0} \sin \theta_{0}\\
                             \cos \phi_{0}\\
                    \end{array}
            \right )$ with $\phi_{0} \in (0, \pi)$, $\theta_{0} \in [0, 2\pi)$. For small $\varepsilon > 0$, the periodic
solution $x(t,x_{0})=e^{-X_{0}t}x_{0}$ of the differential equations (\ref{E:pertdif_projected}) obtained for
$\varepsilon= 0$ is continuously deformed when $\varepsilon$ increases to another periodic solution $x^{\varepsilon}(t)$
with period $T(\varepsilon)= \frac{2\pi}{\left | X_{0} \right | }+O(\varepsilon)$ if \[I(\phi_{0})= 0 \mbox{ and }
I^{'}(\phi_{0}) \neq 0.\] Also $x^{0}(t)= x(t,x_{0})$.\\
The stability of the periodic solution $x^{\varepsilon}(t)$ is as follows:
\begin{enumerate}
\item if $I^{'}(\phi_{0}) < 0$, then $x^{\varepsilon}(t)$ is locally asymptotically stable for the differential equations
(\ref{E:pertdif_projected});
\item if $I^{'}(\phi_{0}) > 0$, then $x^{\varepsilon}(t)$ is unstable for the differential equations
(\ref{E:pertdif_projected}).
\end{enumerate}
\end{enumerate}
Using the Poincar\'{e}-Bendixson theorem on a sphere (\cite{HS}), we have a complete picture of the phase portrait for
the differential equations (\ref{E:pertdif_projected}).
\end{theorem}
The matrix $B \in SO(3)$ such that $Bx_{1}^{0}= (0,0,1)^{T}$ is not unique, but it can be shown that
Theorem \ref{thm:FSB_DE_Orbit_E} does not depend on the choice of the matrix $B$. Therefore, in the statement and
proof of Theorem \ref{thm:FSB_DE_Orbit_E}, without loss of generality we may assume that $x_{1}^{0}=(0,0,1)^{T}$.
The proof of Theorem \ref{thm:FSB_DE_Orbit_E} is presented in Section \ref{S:Proof_FSB}.\\
Let $\varepsilon > 0$ be fixed and $\phi_{0} \in (0,\pi)$ such that it corresponds to the $\frac{2\pi}
{\left | X_{0} \right |}$-periodic solution $x(t, x_{0})=e^{-X_{0}t}x_{0}$ of the differential equations
(\ref{E:pertdif_projected}) for $\varepsilon= 0$. Then, the condition $I(\phi_{0})=0$ is necessary for the
persistence of $\frac{2\pi}{\left | X_{0} \right | }$-periodic solution $x(t, x_{0})$ to $T(\varepsilon)$-periodic
solution $x^{\varepsilon}(t)$.

\section{The Analysis of the $SO(2)$-Equivariant Dynamics on the Phase Space $SO(3)$}\label{S:Phase_space}

We translate the results of Section \ref{S:Projection}, which are valid for the orbit space $SO(3)/SO(2)$, to the
phase space $SO(3)$ (Proposition \ref{prop:FSB_Phase_E}) or $\mathbf{Y^{\alpha}}$ (Theorem \ref{thm:results_FSB}).
\begin{proposition}[$\varepsilon > 0$]\label{prop:FSB_Phase_E}
For $\varepsilon > 0$ small we have:
\begin{enumerate}
\item The equilibrium $x^{1}(\varepsilon)$ of the differential equations (\ref{E:pertdif_projected}) gives the
$SO(2)$-orbit of a solution of the differential equations (\ref{E:pertdif}) having the form $A(t,\varepsilon)=
e^{\alpha_{1}(\varepsilon)Q t}A_{\varepsilon}$, where $x^{Q}_{1}= A_{\varepsilon}x^{1}(\varepsilon)$ and
$\alpha_{1}(\varepsilon)= \left | X_{0} \right | +O(\varepsilon)$. The equilibrium $x^{1}(\varepsilon)$ is generically
hyperbolic, therefore the corresponding $SO(2)$-orbit is generically normally hyperbolic;
\item The equilibrium $x^{2}(\varepsilon)$ of the differential equations (\ref{E:pertdif_projected}) gives the
$SO(2)$-orbit of a solution of the differential equations (\ref{E:pertdif}) having the form
$A(t, \varepsilon)= e^{\alpha_{2}(\varepsilon)Qt}A_{\varepsilon}$, where $x^{Q}_{1}= A_{\varepsilon}x^{2}(\varepsilon)$
and $\alpha_{2}(\varepsilon)= -\left | X_{0} \right | +O(\varepsilon)$. The equilibrium $x^{2}(\varepsilon)$ is
generically hyperbolic, therefore the corresponding $SO(2)$-orbit is generically normally hyperbolic;
\item If $I(\phi_{0})= 0$ and $I^{'}(\phi_{0}) \neq 0$, where $I$ and $\phi_{0}$ are defined in
Theorem \ref{thm:FSB_DE_Orbit_E}, then any $T(\varepsilon)$-periodic solution $x^{\varepsilon}(t)$ of
(\ref{E:pertdif_projected}) gives the $SO(2)$-orbit of a solution of the differential equations (\ref{E:pertdif})
having the form $A(t, \varepsilon)= e^{\beta(\varepsilon)Qt}B(t,\varepsilon)$, where $\beta(\varepsilon)=
O(\varepsilon)$ and $B(t,\varepsilon)$ is a $T(\varepsilon)= \frac{2\pi}{\left | X_{0} \right | } +O(\varepsilon)$-
periodic function such that $x^{Q}_{1}= B(0, \varepsilon)x^{\varepsilon}(0)$.
\end{enumerate}
The orbital stability of the above solutions is the same as the stability of the solutions of (\ref{E:pertdif_projected})
from which they were obtained (see \cite{Ch}). Using the Poincar\'{e}-Bendixson theorem on a sphere (\cite{HS})
and \cite{Ch}, we have a complete picture of the phase portrait for the differential equations (\ref{E:pertdif}).
\end{proposition}
The proof of Proposition \ref{prop:FSB_Phase_E} is presented in Section \ref{S:Proof_FSB}.\\
The second main result of this article is the following theorem:
\begin{theorem}[The Effects of Forced Symmetry-Breaking]
\label{thm:results_FSB}
For small $\varepsilon > 0$ we have:
\begin{enumerate}
\item The $SO(2)$-orbit of the rotating wave $\mathbf{\Phi}(t,A_{0}u_{0},0)= e^{Q\left | X_{0} \right |
t}A_{0}u_{0}$ is deformed as $\varepsilon$ increases to the $SO(2)$-orbit of a rotating wave having as the frequency
vector $(\left | X_{0} \right | +O(\varepsilon))\ \overrightarrow Q$, where $\overrightarrow Q=
A_{0}\overrightarrow X_{0}$;
 \item
The $SO(2)$-orbit of the rotating wave $\mathbf{\Phi}(t,A_{0}u_{0},0)= e^{-Q\left | X_{0} \right | t}A_{0}u_{0}$ is
deformed as $\varepsilon$ increases to a rotating wave having as the frequency vector $(-\left | X_{0} \right | +
O(\varepsilon))\overrightarrow Q$;
\item The $SO(2)$-orbit of the rotating waves $\mathbf{\Phi}(t,A_{0}u_{0},0)= e^{X t}A_{0}u_{0}$
having the frequency vectors $\overrightarrow X=\overrightarrow {A_{0}X_{0}A_{0}^{-1}}$ different from
$\left | X_{0} \right | \overrightarrow Q$ and $-\left |  X_{0} \right | \overrightarrow Q$, and satisfying
$I(\phi_{0})= 0$ and $I^{'}(\phi_{0}) \neq 0$, where $I$ and $\phi_{0}$ are defined in Theorem \ref{thm:FSB_DE_Orbit_E},
transform as $\varepsilon$ increases into $SO(2)$-orbits of modulated rotating waves having primary
frequency vectors of the form $O(\varepsilon)\overrightarrow Q$ and the relative period of the form
$T(\varepsilon)=\frac{2\pi}{\left | X_{0} \right |}+O(\varepsilon)$.
\end{enumerate}
The orbital stability of the above solutions of the reaction-diffusion system (\ref{E:FSB_RD}) is related to the
orbital stability of the rotating wave $\mathbf{\Phi}(t,u_{0},0)$ and to the stability of the solutions of
(\ref{E:pertdif_projected}) from which they were obtained as follows:
\begin{itemize}
\item If the rotating wave $\mathbf{\Phi}(t,u_{0},0)$ is orbitally unstable, then the above solutions of the
reaction-diffusion system (\ref{E:FSB_RD}) are orbitally unstable;
\item If the rotating wave $\mathbf{\Phi}(t,u_{0},0)$ is orbitally stable (that is $E^{cu}=T_{u_{0}}(SO(3)u_{0})$),
then the above solutions of the reaction-diffusion system (\ref{E:FSB_RD}) have the same orbital stability as the
stability of the solutions of (\ref{E:pertdif_projected}) from which they were obtained.
\end{itemize}
\end{theorem}
The proof of Theorem \ref{thm:results_FSB} is presented in Section \ref{S:Proof_FSB}. In \cite{Co}, we present some
examples that may explain similar phenomena to those observed for planar spiral waves (boundary drifting, anchoring at
the localized inhomogeneities, repelling from the localized inhomogeneities), and another new phenomenon that appears
due to the curvature of the sphere (anchoring at the opposite site).\\

\section{Proofs of Theorems} \label{S:Proof_FSB}

\begin{proof}[Proof of Proposition \ref{prop:Pert_RE}]
There is nothing to prove for $G_{pert}^{2}(A, q, \varepsilon)$.\\
Let $k(A, q, \varepsilon)= A^{-1}G_{pert}^{1}(A, q, \varepsilon)$ for $A
\in SO(3)$, $q \in V_{*}$ and $\varepsilon \geq 0$ small.\\
Since $G_{pert}^{1}(A, q, \varepsilon) \in T_{A}(SO(3))= A\cdot
so(3)$, then $k(A, q, \varepsilon) \in so(3)$ for any
$A \in SO(3), q \in V_{*}$ and $\varepsilon \geq 0$ small.
Using the $SO(2)$-equivariance of $G_{pert}^{1}(., q, \varepsilon)$,
we get that the function $k(., q, \varepsilon)$ is \\
$SO(2)$-invariant for $\varepsilon \geq 0$ small and $q \in V_{*}$, that is\\
$k(e^{Q\theta}A,q,\varepsilon)=(e^{Q\theta}A)^{-1}G_{pert}^{1}(e^{Q\theta}A,q,\varepsilon)=
A^{-1}e^{-Q\theta}e^{Q\theta}G_{pert}^{1}(A,q,\varepsilon)=
A^{-1}G_{pert}^{1}(A,q,\varepsilon)=k(A,q,\varepsilon)$ for any $A \in
SO(3)$, $q \in V_{*}$, $\theta \in [0,2\pi)$ and $\varepsilon \geq 0$
small.
\end{proof}
\begin{proof}[Proof of Theorem \ref{thm:DE_FSB_RE}]
For $\varepsilon > 0$ small, $M(\varepsilon)$ corresponds to the hyperbolic equilibrium $q_{\varepsilon} \in V_{*}$
of the second differential equation of (\ref{E:pertdifCM}). Therefore, we only consider the $SO(2)$-equivariant dynamics
on $SO(3)$ obtained by substituting $q= q_{\varepsilon}$ in the differential equations (\ref{E:pertdifCM}).
Using Proposition \ref{prop:Pert_RE}, we get
\begin{equation}\label{E:diffaux1}
\dot{A}= A\left [f_{1}(q_{\varepsilon})+\varepsilon k(A,q_{\varepsilon},\varepsilon)\right ].
\end{equation}
Since $q_{\varepsilon}=O(\varepsilon)$ and $f_{1}$ is sufficiently smooth, we get $f_{1}(q_{\varepsilon})= f_{1}(0)+
\varepsilon H(\varepsilon)= X_{0}+\varepsilon H(\varepsilon)$ with $H(\varepsilon) \in so(3)$. The
differential equations (\ref{E:diffaux1}) become
\begin{equation}\label{E:diffaux2}
\dot{A}= A\left [X_{0}+\varepsilon H(\varepsilon)+\varepsilon k(A,q_{\varepsilon},\varepsilon)\right ].
\end{equation}
Let $g(A,\varepsilon)= H(\varepsilon)+ k(A,q_{\varepsilon},\varepsilon)$. We get the differential equations
(\ref{E:pertdif}), where $g(.,\varepsilon)$ is $SO(2)$-invariant, since $k(.,q_{\varepsilon},\varepsilon)$ is
$SO(2)$-invariant.
\end{proof}
\begin{proof}[Proof of Theorem \ref{thm:FSB_PROJ_E}]
It is enough to prove that there exists a sufficiently smooth function $g^{S} \colon \mathbf{S^{2}} \times
[0,\varepsilon_{0}) \rightarrow so(3)$ for $\varepsilon_{0}$ small such that
\begin{equation}\label{E:ProjectedPert}
\widetilde{g}(C, \varepsilon)= g^{S}(Cx^{1}_{Q}, \varepsilon) \mbox { for any } C\in SO(3) \mbox { and }
\varepsilon \geq 0 \mbox{ small }.
\end{equation}
Since $\mathbf{S^{2}}$ is compact, it is enough to prove (\ref{E:ProjectedPert}) locally on $\mathbf{S^{2}}$. We may then
use the partition of unity to get a sufficiently smooth function $g^{S}$ globally defined on $\mathbf{S^{2}}$.\\
Let $C_{0} \in SO(3)$ be fixed, but arbitrary. There exists a smooth local section $\pi_{s}$ of $SO(3)/SO(2)$ in $SO(3)$
near $C_{0} \cdot SO(2)$ (see Theorem \ref{thm:Orbitspace_Manifold}). We have:
\begin{equation}\label{E:AuxFSB}
\begin{split}
\pi \circ \pi_{s} & = id  \\ & \Rightarrow \pi(\pi_{s}(C\cdot SO(2)))= C\cdot SO(2)\\
& \Rightarrow \pi_{s}(C\cdot SO(2))SO(2)= C\cdot SO(2) \\
& \Rightarrow C=\pi_{s}(C\cdot SO(2))\cdot C^{*}(C) \mbox{ for } C^{*}(C) \in SO(2) \mbox{ near } C_{0} \cdot SO(2).
\end{split}
\end{equation}
Let us define the function
\begin{equation}\label{E:SphereFunction}
g^{S}= \widetilde{g} \circ\pi_{s}\circ\beta^{-1} \mbox{ near } C_{0}x_{Q}^{1}.
\end{equation}
Then $g^{S}$ is a sufficiently smooth function defined locally near $C_{0}x_{Q}^{1}$ on $\mathbf{S^{2}}$. Also,
\begin{equation}
g^{S}(Cx_{Q}^{1}, \varepsilon)= \widetilde{g}(\pi_{s}(C\cdot SO(2)),\varepsilon)=
\widetilde{g}(\pi_{s}(C\cdot SO(2))\cdot C^{*}(C), \varepsilon)= \widetilde{g}(C, \varepsilon)
\end{equation}
for $C$ near $C_{0}$ and $\varepsilon \geq 0$ small, where for the second equality we use the invariance of
$\widetilde{g}(., \varepsilon)$ under the restriction of the action $\widetilde{\theta}$ to $SO(2)$ and for the
third equality we use $(\ref{E:AuxFSB})$.\\
The function $g^{S}$ is well-defined since the function $\widetilde{g}$ is $SO(2)$-invariant under the action
$\widetilde{\theta}$ restricted to $SO(2)$; that is, for any $\theta \in [0,2\pi)$, $C \in SO(3)$ and
$\varepsilon \geq 0$ small, we have $g^{S}(Ce^{Q\theta}x_{Q}^{1},\varepsilon)=\widetilde{g}(Ce^{Q\theta},\varepsilon)=
\widetilde{g}(C,\varepsilon)=g^{S}(Cx_{Q}^{1},\varepsilon)$.\\
Also, taking into account that the action of $SO(3)$ on $\mathbf{S^{2}}$ is transitive, we can look at the function
$g^{S}$ as being defined on $\mathbf{S^{2}}$.\\
We check that the vector field given by the differential equations (\ref{E:pertdif_projected}) is associated
to the semiflow $\widetilde{\mathbf{\Psi}_{1}}$. We have $\widetilde{\mathbf{\Psi}_{1}}(t,Cx_{Q}^{1},\varepsilon)=
(\beta \circ \pi)(\mathbf{\Psi}(t,C,\varepsilon))= \mathbf{\Psi}(t,C,\varepsilon)x_{Q}^{1}$. Therefore,
\begin{equation}
\begin{split}
\frac{\partial \widetilde{\mathbf{\Psi}_{1}}}{\partial t}(t,Cx_{Q}^{1},\varepsilon) &=
\frac{\partial}{\partial t}(\mathbf{\Psi}(t,C,\varepsilon)x_{Q}^{1})=
\frac{\partial\mathbf{\Psi}}{\partial t}(t,C,\varepsilon)x_{Q}^{1}\\ &=
-\left [X_{0}+\varepsilon \widetilde{g}(\mathbf{\Psi}(t,C,\varepsilon),\varepsilon)\right ]
\mathbf{\Psi}(t,C,\varepsilon) x_{Q}^{1} \\ & =
 -\left [X_{0}+\varepsilon g^{S}(\mathbf{\Psi}(t,C,\varepsilon)x_{Q}^{1},\varepsilon) \right ]
 \mathbf{\Psi}(t,C,\varepsilon)x_{Q}^{1}\\ & =
-\left [X_{0}+\varepsilon g^{S}(\widetilde{\mathbf{\Psi}_{1}}(t,C x_{Q}^{1},\varepsilon),\varepsilon) \right ]
\widetilde{\mathbf{\Psi}_{1}}(t,C x_{Q}^{1},\varepsilon)
\end{split}
\end{equation}
for any $C \in SO(3)$, $t \in [0,\infty)$ and $\varepsilon \geq 0$ small.
\end{proof}
\begin{proof}[Proof of Proposition \ref{prop:FSB_Orbit_Phase_E}]
Conclusions (1) and (2) are obvious. The third conclusion results from the definition of the flows
$\mathbf{\Psi}$ and $\widetilde{\mathbf{\Psi}_{1}}$, as well as $\mathbf{\Psi}(t,C,\varepsilon)=
[\mathbf{\Phi}(t,C^{-1},\varepsilon)]^{-1}$.
\end{proof}
\begin{proof}[Proof of Theorem \ref{thm:RotatingWaves}]
Conclusions (1) and (2) result from Theorem \ref{thm:CMR_FSB_RE} and Proposition \ref{prop:FSB_Orbit_Phase_E}.\\
Using Theorem \ref{thm:CMR_FSB_RE}, we get that the rotating waves for the reaction-diffusion system (\ref{E:FSB_RD})
for $\varepsilon= 0$ are of the form $\mathbf{\Phi}(t,A_{0}u_{0},0)=A_{0}e^{X_{0}t}u_{0}$ for $A_{0} \in SO(3)$ and
$t \in [0,\infty)$.\\
Therefore, the primary frequency vector of $\mathbf{\Phi}(t,A_{0}u_{0},0)$ is $\overrightarrow{A_{0}X_{0}A_{0}^{-1}}$.
By Proposition \ref{prop:FSB_Orbit_Phase_E}, the rotating wave $\mathbf{\Phi}(t,A_{0}u_{0},0)$ projects onto
$x_{0} \in \mathbf{S}^{2}$ such that $x_{1}^{Q}=A_{0}x_{0}$.\\
Let $X= \frac{1}{\left | X_{0} \right | }A_{0}X_{0}A_{0}^{-1}$. Then $x_{1}^{0}=\frac{1}{\left | X_{0} \right |}
\overrightarrow X_{0}$, we get $A_{0}x_{1}^{0}=\overrightarrow{X}$. Also, $(x_{1}^{Q})^{T}=x_{0}^{T}A_{0}^{T}$.
Therefore, $\cos(\angle(\vec X, \vec Q))=(x_{1}^{Q})^{T}\overrightarrow{X}= x_{0}^{T}A_{0}^{T}A_{0}x_{1}^{0}=
x_{0}^{T}x_{1}^{0}= \cos(\angle(\vec{x_{0}},\vec{X_{0}}))$.\\
The third conclusion results from the first conclusion.
\end{proof}
\begin{proof}[Proof of Theorem \ref{thm:FSB_DE_Orbit_E}]
For $\varepsilon = 0$, the differential equation (\ref{E:pertdif_projected}) has two equilibria, $x^{0}_{1}$ and
$x^{0}_{2}$. Using the Implicit Function Theorem , we will prove the persistence of these two equilibria for small
$\varepsilon > 0$.\\
We have defined $B \in SO(3)$ by $Bx_{1}^{0}= (0,0,1)^{T}$. Then $Bx_{2}^{0}= (0,0,-1)^{T}$. Let $s= Bx$. Then,
\begin{equation}\label{E:EqDiffNew1}
\dot s= B\dot x= -\left [BX_{0}B^{-1}+\varepsilon Bg^{S}(B^{-1}s, \varepsilon)B^{-1} \right ]s,
\end{equation}
and, if we write $X_{00}=BX_{0}B^{-1}$ and $g^{SS}(s,\varepsilon)= Bg^{S}(B^{-1}s, \varepsilon)B^{-1}$, then we get
\begin{equation}\label{E:EqDiffNew2}
\dot s= -\left [X_{00}+\varepsilon g^{SS}(s, \varepsilon)\right ]s,
\end{equation}
Since $B\overrightarrow {X_{0}} =\left | X_{0} \right | (0,0,1)^{T}$, we have
\[X_{00}= BX_{0}B^{-1}= \left | X_{0} \right | \left (
                                                     \begin{array}{ccc}
                                                      0 & -1 & 0 \\
                                                      1 &  0 & 0 \\
                                                      0 &  0 & 0  \\
                                                     \end{array}
                                                \right ).\]
Let us define the function
\begin{equation}\label{E:FunctionImplicit}
F(s, \varepsilon)= -\left [X_{00}+\varepsilon g^{SS}(s, \varepsilon)\right]s, \mbox { where }
F \colon \mathbf{S^{2}} \times [0,\varepsilon_{0}) \rightarrow T\mathbf{S^{2}} \mbox{ for }\varepsilon_{0} >0
\mbox{ small }.
\end{equation}
It is clear that $F((0,0,1), 0)= 0$. We will compute $(D_{s}F)_{((0,0,1),0)} \colon T_{(0,0,1)}\mathbf{S^{2}}
\rightarrow T\mathbf{S^{2}}\simeq \mathbb{R}^{2}$.\\
Let $s= (x_{1}, x_{2}, x_{3}) \in \mathbf{S^{2}}$ and recall that
\[g^{SS}(s,\varepsilon)= g^{SS}(x_{1}, x_{2}, x_{3},\varepsilon)=
\left ( \begin{array}{ccc}
0 & -G_{1}(x_{1}, x_{2}, x_{3},\varepsilon) & G_{2}(x_{1}, x_{2}, x_{3},\varepsilon)\\
G_{1}(x_{1}, x_{2}, x_{3},\varepsilon) & 0  & -G_{3}(x_{1}, x_{2}, x_{3},\varepsilon)\\
-G_{2}(x_{1}, x_{2}, x_{3},\varepsilon) & G_{3}(x_{1}, x_{2}, x_{3},\varepsilon) & 0\\
\end{array} \right ). \]
For $(x_{1}, x_{2}, x_{3}) \in \mathbf{S^{2}}$ near $(0,0,1)$, one has $x_{3}= \sqrt{1-x_{1}^{2}-x_{2}^{2}} > 0$
for $\left \| (x_{1}, x_{2}) \right \| \geq 0$ small enough. If we substitute this expression of $x_{3}$ into the
differential equations (\ref{E:EqDiffNew2}), we have:
\begin{multline}\label{E:Diff_aux1}
\left ( \begin{array}{c}
          \dot{x_{1}}\\
          \dot{x_{2}}\\
          \frac{d }{d t} \sqrt{1-x_{1}^{2}-x_{2}^{2}}\\
         \end{array}
\right )= -\left[\left ( \begin{array}{ccc}
                  0 & -\left | X_{0} \right | & 0\\
                  \left | X_{0} \right | & 0 & 0\\
                  0& 0& 0\\
                  \end{array}
           \right )\right .\\ \left . +
\varepsilon \left ( \begin{array}{ccc}
0 & -G_{1}(x_{1}, x_{2}, \sqrt{1-x_{1}^{2}-x_{2}^{2}} ,\varepsilon) &
G_{2}(x_{1}, x_{2}, \sqrt{1-x_{1}^{2}-x_{2}^{2}},\varepsilon)\\
G_{1}(x_{1}, x_{2}, \sqrt{1-x_{1}^{2}-x_{2}^{2}},\varepsilon) & 0  &
-G_{3}(x_{1}, x_{2}, \sqrt{1-x_{1}^{2}-x_{2}^{2}},\varepsilon)\\
-G_{2}(x_{1}, x_{2}, \sqrt{1-x_{1}^{2}-x_{2}^{2}},\varepsilon) &
G_{3}(x_{1}, x_{2}, \sqrt{1-x_{1}^{2}-x_{2}^{2}},\varepsilon) & 0\\
\end{array} \right ) \right] \\ \cdot \left ( \begin{array}{c} x_{1}\\ x_{2}\\ \sqrt{1-x_{1}^{2}-x_{2}^{2}}\\ \end{array}
\right )
\end{multline}
or
\begin{equation}\label{E:Diff_aux2}
\begin{array}{lll}
\dot{x_{1}} &=& \left | X_{0} \right | x_{2}+\varepsilon G_{1}(x,\varepsilon)x_{2}-
\varepsilon G_{2}(x,\varepsilon)\sqrt{1-x_{1}^{2}-x_{2}^{2}},\\
\dot{x_{2}} &=& -\left | X_{0} \right | x_{1}-\varepsilon G_{1}(x,\varepsilon)x_{1}+
\varepsilon G_{3}(x,\varepsilon) \sqrt{1-x_{1}^{2}-x_{2}^{2}},\\
\frac{d}{d t }\sqrt{1-x_{1}^{2}-x_{2}^{2}} &=& \varepsilon G_{2}(x,\varepsilon)x_{1}-
\varepsilon G_{3}(x,\varepsilon)x_{2},\\
\end{array},
\end{equation}
where $x= (x_{1}, x_{2}, \sqrt{1-x_{1}^{2}-x_{2}^{2}})$.\\
If the first two equations in the system (\ref{E:Diff_aux2}) are satisfied, then the third one is immediately satisfied.
Therefore, we get the system
\begin{equation}\label{E:EqDiffNew3}
\begin{array}{cc}
\dot{x_{1}} &= H_{1}(x_{1}, x_{2}, \varepsilon),\\
\dot{x_{2}} &= H_{2}(x_{1}, x_{2}, \varepsilon),\\
\end{array}
\end{equation}
where
\begin{equation}\label{E:AuxH1}
\begin{split}
H_{1}(x_{1}, x_{2}, \varepsilon) &= \left | X_{0} \right |x_{2}+
\varepsilon \left [G_{1}(x_{1}, x_{2}, \sqrt{1-x_{1}^{2}-x_{2}^{2}}, \varepsilon)x_{2} \right . \\  &
\left . -G_{2}(x_{1}, x_{2}, \sqrt{1-x_{1}^{2}-x_{2}^{2}}, \varepsilon) \sqrt{1-x_{1}^{2}-x_{2}^{2}}\right ],
\end{split}
\end{equation}
and
\begin{equation}\label{E:AuxH2}
\begin{split}
H_{2}(x_{1}, x_{2}, \varepsilon) &= -\left | X_{0} \right |x_{1}-\varepsilon \left [G_{1}(x_{1}, x_{2},
\sqrt{1-x_{1}^{2}-x_{2}^{2}}, \varepsilon)x_{1} \right .\\ & \left . -G_{3}(x_{1}, x_{2}, \sqrt{1-x_{1}^{2}-x_{2}^{2}},
\varepsilon) \sqrt{1-x_{1}^{2}-x_{2}^{2}} \right ].
\end{split}
\end{equation}
Let $H(x_{1}, x_{2},\varepsilon)= (H_{1}(x_{1}, x_{2},\varepsilon), H_{2}(x_{1}, x_{2},\varepsilon))$. We have
$H(0,0,0)= (0,0)$ and the linearization of $H(x_{1}, x_{2}, \varepsilon)$ about $(0,0)$ at $\varepsilon=0$ is
\[D_{(x_{1}, x_{2})}H(0,0,0)= \left (
                                  \begin{array}{cc}
                                  0  & \left | X_{0} \right |\\
                                  -\left | X_{0} \right | & 0 \\
                                  \end{array}
                              \right ),
\]
which is an invertible matrix. Using the implicit function theorem there exists a sufficiently smooth branch
$(x_{1}^{1}(\varepsilon), x_{2}^{1}(\varepsilon))$ near $(0,0)$ such that $(x_{1}^{1}(0), x_{2}^{1}(0))= (0,0)$,
and $H_{1}(x_{1}^{1}(\varepsilon), x_{2}^{1}(\varepsilon), \varepsilon)=0$ and $H_{2}(x_{1}^{1}(\varepsilon),
x_{2}^{1}(\varepsilon), \varepsilon)=0$ for $\varepsilon \geq 0$ small.\\
We have
\begin{equation}\label{E:Implicit_diff}
x_{1}^{1}(\varepsilon)= x_{11}\varepsilon+O(\varepsilon^{2}) \mbox { and }
x_{2}^{1}(\varepsilon)= x_{22}\varepsilon+O(\varepsilon^{2}), \, \mbox{ for
small } \varepsilon \geq 0.
\end{equation}
By using implicit differentiation we get
\begin{equation}\label{E:Implicit1}
\begin{array}{ccc}
\frac{\partial{H_{1}}}{\partial{x_{1}}}\frac{dx_{1}}{d \varepsilon}+
\frac{\partial{H_{1}}}{\partial{x_{2}}}\frac{dx_{2}}{d \varepsilon} &=&-\frac{\partial{H_{1}}}{\partial{\varepsilon}},\\
\frac{\partial{H_{2}}}{\partial{x_{1}}}\frac{dx_{1}}{d \varepsilon}+
\frac{\partial{H_{2}}}{\partial{x_{2}}}\frac{dx_{2}}{d \varepsilon} &=&-\frac{\partial{H_{2}}}{\partial{\varepsilon}}.\\
\end{array}
\end{equation}
Generically, we have
\begin{equation}\label{E:generic1}
G_{2}(0,0,1,0) \neq 0 \mbox { and } G_{3}(0,0,1,0) \neq 0.
\end{equation}
From (\ref{E:Implicit1}) it follows that at $\varepsilon= 0$ we have
\begin{equation}
\frac{dx_{1}}{d \varepsilon} = -\frac{\left | \begin{array}{cc}
                                        \frac{\partial H_{1}}{\partial \varepsilon} &
                                        \frac{\partial H_{1}}{\partial x_{2}}\\
                                        \frac{\partial H_{2}}{\partial \varepsilon} &
                                        \frac{\partial H_{2}}{\partial x_{2}}\\
                                               \end{array}
                                       \right |}{\left | \begin{array}{cc}
                                        \frac{\partial H_{1}}{\partial x_{1}} &
                                        \frac{\partial H_{1}}{\partial x_{2}}\\
                                        \frac{\partial H_{2}}{\partial x_{1}} &
                                        \frac{\partial H_{2}}{\partial x_{2}}\\
                                                          \end{array}
                                                  \right | }
\end{equation}
and
\begin{equation}
\frac{dx_{2}}{d \varepsilon}= -\frac{\left | \begin{array}{cc}
                                        \frac{\partial H_{1}}{\partial x_{1}} &
                                        \frac{\partial H_{1}}{\partial \varepsilon}\\
                                        \frac{\partial H_{2}}{\partial x_{1}} &
                                        \frac{\partial H_{2}}{\partial \varepsilon}\\
                                         \end{array}
                                      \right |}{\left | \begin{array}{cc}
                                        \frac{\partial H_{1}}{\partial x_{1}} &
                                        \frac{\partial H_{1}}{\partial x_{2}}\\
                                        \frac{\partial H_{2}}{\partial x_{1}} &
                                        \frac{\partial H_{2}}{\partial x_{2}}\\
                                                       \end{array}
                                                \right | }.
\end{equation}
Since $\frac{\partial H_{1}}{\partial x_{1}}(0,0,0)=\frac{\partial H_{2}}{\partial x_{2}}(0,0,0)=0$,
$\frac{\partial H_{1}}{\partial x_{2}}(0,0,0)=-\frac{\partial H_{2}}{\partial x_{1}}(0,0,0)=
\left | X_{0} \right | $, $\frac{\partial H_{1}}{\partial \varepsilon}(0,0,0)=-G_{2}(0,0,1,0)$ and
$\frac{\partial H_{2}}{\partial \varepsilon}(0,0,0)=G_{3}(0,0,1,0)$, it follows that
\begin{equation}\label{E:Implicit2}
\begin{array}{ccc}
x_{11} &=& \frac{1}{\left | X_{0} \right | }G_{3}(0,0,1,0),\\
x_{22} &=& \frac{1}{\left | X_{0} \right | }G_{2}(0,0,1,0).\\
\end{array}
\end{equation}
So, for $\varepsilon \geq 0$ small, we get
\begin{equation}\label{E:Impliciteq}
x_{1}^{1}(\varepsilon)= \frac{1}{\left | X_{0} \right |}G_{3}(0,0,1,0)\varepsilon+
O(\varepsilon^{2}), x_{2}^{1}(\varepsilon)=\frac{1}{\left | X_{0} \right |}G_{2}(0,0,1,0)\varepsilon+O(\varepsilon^{2}).
\end{equation}
Therefore, there exists a sufficiently smooth branch
$s^{1}(\varepsilon)= (x_{1}^{1}(\varepsilon), x_{2}^{1}(\varepsilon), x_{3}^{1}(\varepsilon))$ near
$(0,0,1)$ of equilibria of the differential equations (\ref{E:EqDiffNew2}), such that \\
$s^{1}(0)= (x_{1}^{1}(0), x_{2}^{1}(0), x_{3}^{1}(0))= (0,0,1)$ and $F(s^{1}(\varepsilon), \varepsilon)= (0,0,0)$ for
$\varepsilon \geq 0$ small.\\
For $\varepsilon \geq 0$ small, there is a sufficiently smooth branch $x^{1}(\varepsilon) \in \mathbf{S^{2}}$ of
equilibria for the differential equations (\ref{E:pertdif_projected}) such that $x^{1}(0)= x_{1}^{0}$.\\
In the same way, for $\varepsilon \geq 0$ small there exists a sufficiently smooth branch $x^{2}(\varepsilon)$ of
equilibria for the differential equations (\ref{E:pertdif_projected}) such that $x^{2}(0)= x_{2}^{0}$.\\
The stability of the equilibrium $x^{1}(\varepsilon)$ is the same as the stability of the equilibrium
$(x^{1}_{1}(\varepsilon), x^{1}_{2}(\varepsilon))$ for $\varepsilon \geq 0$ small, and the stability of the equilibrium
$x^{2}(\varepsilon)$ is the same as the stability of the equilibrium $(x^{2}_{1}(\varepsilon), x^{2}_{2}(\varepsilon))$
for $\varepsilon \geq 0$ small.\\
Let us now obtain the stability of the equilibrium $(x^{1}_{1}(\varepsilon), x^{1}_{2}(\varepsilon))$.\\
For $\varepsilon \geq 0$ small, the linearization of $H(x_{1}, x_{2}, \varepsilon)$ around $(x_{1}^{1}(\varepsilon),
x_{2}^{1}(\varepsilon))$ is given by
\begin{equation}\label{E:Linearization}
D_{(x_{1}, x_{2})}H(x_{1}^{1}(\varepsilon), x_{2}^{1}(\varepsilon), \varepsilon)=
     \left ( \begin{array}{cc}
           a_{11}(\varepsilon) & \left | X_{0} \right | +O(\varepsilon)\\
           -\left | X_{0} \right | +O(\varepsilon) & a_{22}(\varepsilon)\\
            \end{array}
     \right ),
\end{equation}
where
\begin{equation}\label{E:first}
a_{11}(\varepsilon)= \frac{\partial H_{1}}{\partial x_{1}}(x_{1}^{1}(\varepsilon), x_{2}^{1}(\varepsilon),\varepsilon)=
a_{11}\varepsilon+O(\varepsilon^{2}) \mbox{ since } a_{11}(0)=0
\end{equation}
and
\begin{equation}\label{E:second}
a_{22}(\varepsilon)= \frac{\partial H_{2}}{\partial x_{2}}(x_{1}^{1}(\varepsilon), x_{2}^{1}(\varepsilon),\varepsilon)=
a_{22}\varepsilon+O(\varepsilon^{2}) \mbox{ since } a_{22}(0)=0.
\end{equation}
The eigenvalues of the linearization of $H(x_{1}, x_{2}, \varepsilon)$ around
$(x_{1}^{1}(\varepsilon), x_{2}^{1}(\varepsilon))$ satisfies
\begin{equation}\label{E:eigenvalues1}
\left |
\begin{array}{cc}
a_{11}(\varepsilon)-\lambda & \left | X_{0} \right | +O(\varepsilon)\\ -\left | X_{0} \right | +O(\varepsilon) &
a_{22}(\varepsilon)-\lambda\\
\end{array}
\right |= 0.
\end{equation}
Then, it follows that
$(a_{11}(\varepsilon)-\lambda)(a_{22}(\varepsilon)-\lambda)-(\left | X_{0} \right | +O(\varepsilon))
(-\left | X_{0} \right | +O(\varepsilon)=0$ or \\
$\lambda^{2}-(a_{11}(\varepsilon)+a_{22}(\varepsilon))\lambda+\left | X_{0} \right | ^{2}+O(\varepsilon)=0$.
We get the eigenvalues
\begin{equation}\label{E:eigenvalues2}
\lambda_{1,2}(\varepsilon)= \frac{a_{11}(\varepsilon)+a_{22}(\varepsilon)}{2} \pm i(\left | X_{0} \right | +
O(\varepsilon))=  \left [ \frac{a_{11}+a_{22}}{2}\varepsilon+O(\varepsilon^{2})\right ] \pm i(\left | X_{0} \right | +
O(\varepsilon)).
\end{equation}
Generically, $a_{11}+a_{22} \neq 0$. In fact, $(x_{1}^{1}(\varepsilon), x_{2}^{1}(\varepsilon))$ is generically a
hyperbolic equilibrium for the differential equations (\ref{E:pertdif_projected}).
If $a_{11}+a_{22} < 0$, then $(x_{1}^{1}(\varepsilon), x_{2}^{1}(\varepsilon))$ is an asymptotically stable equilibrium
for the differential equations (\ref{E:pertdif_projected}), and if $a_{11}+a_{22} > 0$, then $(x_{1}^{1}(\varepsilon),
x_{2}^{1}(\varepsilon))$ is an unstable equilibrium for the differential equations (\ref{E:pertdif_projected}).\\
We express $a_{11}$ and $a_{22}$ in terms of $G_{1}((0,0,1),0)$, $G_{2}((0,0,1),0)$ and $G_{3}((0,0,1),0)$ and their
partial derivatives at $(0,0,1)$.\\
The Taylor expansion of $H_{1}(x_{1}, x_{2}, \varepsilon)$ about $(0,0,0)$ up to terms of order 2 (h.o.t. contains all
the terms of order $\geq 3$ in $x_{1},x_{2},\varepsilon$) is:
\begin{equation}\label{E:TaylorH1}
\begin{split}
H_{1}(x_{1}, x_{2}, \varepsilon)&=
a_{0}^{1}+a_{1}^{1}x_{1}+a_{2}^{1}x_{2}+a_{3}^{1}x_{1}^{2}+a_{4}^{1}x_{2}^{2}+ a_{5}^{1}x_{1}x_{2}\\ & +
a_{6}^{1}\varepsilon+a_{7}^{1}x_{1}\varepsilon+a_{8}^{1}x_{2}\varepsilon+a_{9}^{1}\varepsilon^{2}+h.o.t.
\end{split}
\end{equation}
Since $a_{0}^{1}= 0$, $a_{1}^{1}= 0$, $a_{3}^{1}= 0$, $a_{4}^{1}=0$ and $a_{5}^{1}= 0$, we get
\begin{equation}\label{E:TaylorH11}
H_{1}(x_{1}, x_{2}, \varepsilon)= a_{2}^{1}x_{2}+
a_{6}^{1}\varepsilon+a_{7}^{1}x_{1}\varepsilon+a_{8}^{1}x_{2}\varepsilon+a_{9}^{1}\varepsilon^{2}+h.o.t.
\end{equation}
The Taylor expansion of $H_{2}(x_{1}, x_{2}, \varepsilon)$ about $(0,0,0)$ up to terms of order 2 is:
\begin{equation}\label{E:TaylorH2}
\begin{split}
H_{2}(x_{1}, x_{2}, \varepsilon) &= a_{0}^{2}+a_{1}^{2}x_{1}+a_{2}^{2}x_{2}+a_{3}^{2}x_{1}^{2}+a_{4}^{2}x_{2}^{2}+
a_{5}^{2}x_{1}x_{2}\\ & + a_{6}^{2}\varepsilon+a_{7}^{2}x_{1}\varepsilon+a_{8}^{2}x_{2}\varepsilon+
a_{9}^{2}\varepsilon^{2}+h.o.t.
\end{split}
\end{equation}
Since $a_{0}^{2}= 0$, $a_{2}^{2}= 0$, $a_{3}^{2}= 0$, $a_{4}^{2}= 0$ and $a_{5}^{2}= 0$, we get
\begin{equation}\label{E:TaylorH22}
H_{2}(x_{1}, x_{2}, \varepsilon)= a_{1}^{2}x_{1}+ a_{6}^{2}\varepsilon+a_{7}^{2}x_{1}\varepsilon+
a_{8}^{2}x_{2}\varepsilon+a_{9}^{2}\varepsilon^{2}+h.o.t.
\end{equation}
We have
\begin{equation}\label{E:firstcoeff}
\frac{\partial{H_{1}}}{\partial{x_{1}}}(x_{1}^{1}(\varepsilon), x_{2}^{1}(\varepsilon),\varepsilon)=
a_{7}^{1}\varepsilon+o(\varepsilon^{2})
\end{equation}
and
\begin{equation}\label{E:secondcoeff}
\frac{\partial{H_{2}}}{\partial{x_{2}}}(x_{1}^{1}(\varepsilon), x_{2}^{1}(\varepsilon),\varepsilon)=
a_{8}^{2}\varepsilon+o(\varepsilon^{2}).
\end{equation}
We have that
\begin{equation}\label{E:Stability1}
a_{11}= a_{7}^{1}= \frac{\partial^{2}{H_{1}}}{\partial{x_{1}}\partial{\varepsilon}}(0,0,0)=
-\frac{\partial{G_{2}}}{\partial{x_{1}}}(0,0,1,0),
\end{equation}
by the definition of  $H_{1}$ and
\begin{equation}\label{E:Stability2}
a_{22}= a_{8}^{2}= \frac{\partial^{2}{H_{2}}}{\partial{x_{2}}\partial{\varepsilon}}(0,0,0)=
\frac{\partial{G_{3}}}{\partial{x_{2}}}(0,0,1,0),
\end{equation}
by the definition of $H_{2}$.\\
Generically, we have $a_{7}^{1}+a_{8}^{2} \neq 0$. If $a_{7}^{1}+a_{8}^{2} < 0$, then $(x_{1}^{1}(\varepsilon),
x_{2}^{1}(\varepsilon))$ is an asymptotically stable equilibrium for the differential equations (\ref{E:pertdif_projected}),
and if $a_{7}^{1}+a_{8}^{2} > 0$, then $(x_{1}^{1}(\varepsilon), x_{2}^{1}(\varepsilon))$ is an unstable equilibrium
for the differential equations (\ref{E:pertdif_projected}).\\
In the same way, we can establish the stability of the equilibrium $x^{2}(\varepsilon)$.
This ends the proof of the first conclusion of Theorem \ref{thm:FSB_DE_Orbit_E}.\\
For $\varepsilon= 0$ the differential equation (\ref{E:pertdif_projected}) has $\frac{2\pi} {\left | X_{0} \right |}$-
periodic solutions of the form $x(t, x_{0})= e^{-X_{0}t}x_{0}$, where  $x_{0} \neq x_{1}^{0}$ and  $x_{0} \neq x_{2}^{0}$.\\
Let us fix a point $x_{0}$ on the unit sphere such that $x_{0} \neq x_{1}^{0}$ and $x_{0} \neq x_{2}^{0}$.\\
We will prove the persistence of the periodic solution $x(t, x_{0})$ if some conditions are satisfied.\\
Making the same change of variable as before, we get the system (\ref{E:EqDiffNew2}) and the periodic solution
$x(t, x_{0})$ becomes $s(t, s_{0})= e^{-X_{00}t}s_{0}$, where $s_{0}= Bx_{0}$ and $s_{0} \neq (0,0,1)$,
$s_{0} \neq (0,0,-1)$.\\
For $s \neq (0,0,1)$ and $s \neq (0,0,-1)$ we use the spherical polar coordinates $(\phi, \theta)$ on $\mathbf{S^{2}}$.\\
Let $s= \left ( \begin{array}{c}
                 \sin \phi \cos \theta \\
                 \sin \phi \sin \theta\\
                 \cos \phi\\
                 \end{array}
         \right )$, where $\phi \in (0,\pi)$ and $\theta \in [0, 2\pi)$. Recall that
$s_{0}= \left ( \begin{array}{c}
                 \sin \phi_{0}\cos \theta_{0} \\
                 \sin \phi_{0} \sin \theta_{0}\\
                 \cos \phi_{0}\\
                 \end{array}
            \right )$\\
and \[a(\phi, \theta, \varepsilon)= G_{1}(\sin \phi \cos \theta, \sin \phi \sin \theta, \cos \phi,\varepsilon),\]
\[b(\phi, \theta, \varepsilon)= G_{2}(\sin \phi \cos \theta, \sin \phi \sin \theta, \cos \phi,\varepsilon),\]
\[c(\phi, \theta, \varepsilon)= G_{3}(\sin \phi \cos \theta, \sin \phi \sin \theta, \cos \phi,\varepsilon)\]
for $\theta \in [0, 2\pi)$, $\phi \in [0, \pi]$, $\varepsilon \geq 0$ small and such that
\[a(0, \theta, \varepsilon)= G_{1}(0,0,1,\varepsilon), a(\pi, \theta, \varepsilon)=
 G_{1}(0,0,-1,\varepsilon) \mbox {for any } \theta \in [0, 2\pi),\]
\[b(0, \theta, \varepsilon)= G_{2}(0,0,1,\varepsilon), b(\pi, \theta, \varepsilon)=
G_{2}(0,0,-1,\varepsilon) \mbox {for any } \theta \in [0, 2\pi),\]
\[c(0, \theta, \varepsilon)= G_{3}(0,0,1,\varepsilon), c(\pi, \theta, \varepsilon)=
G_{13}(0,0,-1,\varepsilon) \mbox {for any } \theta \in [0, 2\pi).\] The system (\ref{E:EqDiffNew2}) in
$(\phi, \theta)$ coordinates yields:
\begin{equation*}
\begin{split}
\left ( \begin{array}{c}
          \frac{d}{d t}(\sin \phi \cos \theta )\\
          \frac{d}{d t}(\sin \phi \sin \theta)\\
          \frac{d}{d t}(\cos \phi)\\
         \end{array}
\right ) &= -\left [ \left ( \begin{array}{ccc}
                  0 & -\left | X_{0} \right | & 0\\
                  \left | X_{0} \right | & 0 & 0\\
                  0& 0& 0\\
                    \end{array}
            \right )\right . \\ & \left . +\varepsilon \left ( \begin{array}{ccc}
                                         0 & -a(\phi, \theta, \varepsilon) & b(\phi, \theta, \varepsilon)\\
                                         a(\phi, \theta, \varepsilon)& 0 & -c(\phi, \theta, \varepsilon)\\
                                         -b(\phi, \theta, \varepsilon)& c(\phi, \theta, \varepsilon)& 0\\
                                                            \end{array}
                                                   \right )\right ]
\left ( \begin{array}{c}
          \sin \phi \cos \theta\\
          \sin \phi \sin \theta\\
          \cos \phi\\
         \end{array}
\right )
\end{split}
\end{equation*}
$\Longrightarrow$
\begin{multline*}
\left ( \begin{array}{c}
          \dot{\phi}\cos \phi \cos \theta-\dot{\theta}\sin \phi \sin \theta\\
          \dot{\phi}\cos \phi \sin \theta+\dot{\theta}\sin \phi \cos \theta\\
          -\dot{\phi}\cos \phi\\
         \end{array}
\right )\\=-\left ( \begin{array}{ccc}
                  0 & -\left | X_{0} \right | -\varepsilon a(\phi, \theta, \varepsilon)&
                  \varepsilon b(\phi, \theta, \varepsilon)\\
                  \left | X_{0} \right | + \varepsilon a(\phi, \theta, \varepsilon)& 0 &
                  -\varepsilon c(\phi, \theta, \varepsilon)\\
                  -\varepsilon b(\phi, \theta, \varepsilon)& \varepsilon c(\phi, \theta, \varepsilon)& 0\\
                    \end{array}
             \right )
 \left ( \begin{array}{c}
          \sin \phi \cos \theta\\
          \sin \phi \sin \theta\\
          \cos \phi\\
               \end{array}
       \right )
\end{multline*}
$\Longrightarrow$
\begin{equation}\label{E:EqDiffNew13}
\begin{array}{rll}
\dot{\phi}\cos \phi \cos \theta-\dot{\theta}\sin \phi \sin \theta
&=& \left [\left | X_{0} \right | +\varepsilon a(\phi, \theta,
\varepsilon)\right ]\sin \phi \sin \theta -\varepsilon b(\phi, \theta, \varepsilon)\cos \phi,\\
\dot{\phi}\cos \phi \sin \theta+\dot{\theta}\sin \phi \cos \theta
&=& -\left [\left | X_{0} \right | +\varepsilon a(\phi, \theta,
\varepsilon)\right ]\sin \phi \cos \theta + \varepsilon c(\phi, \theta, \varepsilon)\cos \phi,\\
-\dot{\phi}\sin \phi &=& \varepsilon b(\phi, \theta, \varepsilon)\sin
\phi \cos \theta -\varepsilon c(\phi, \theta, \varepsilon)
\sin \phi \sin \theta.\\
\end{array}
\end{equation}
We divide the third equation in (\ref{E:EqDiffNew13}) by $\sin \phi \neq 0$, where $\phi \in (0,\pi)$ to get
\begin{equation}\label{E:EqDiffNew14}
\dot{\phi}= -\varepsilon b(\phi, \theta, \varepsilon) \cos \theta +\varepsilon c(\phi, \theta, \varepsilon) \sin \theta.
\end{equation}
We substitute $\dot{\phi}$ given by (\ref{E:EqDiffNew14}) into the first two equations of (\ref{E:EqDiffNew13}) to get
\begin{eqnarray}\label{E:EqDiffNew15}
\dot{\theta}\sin \phi \sin \theta = -\left [\left | X_{0} \right | +
\varepsilon a(\phi, \theta, \varepsilon)\right ]\sin \theta \sin \phi \nonumber & +
& \varepsilon c(\phi,\theta,\varepsilon) \cos \theta \sin \theta \cos \phi \nonumber \\+
\varepsilon b(\phi,\theta,\varepsilon)\cos \phi \sin ^{2}\theta, \nonumber \\
\dot{\theta}\sin \phi \cos \theta = -\left [\left | X_{0} \right |
+ \varepsilon a(\phi, \theta, \varepsilon)\right ]\cos \theta \sin \phi
&+& \varepsilon c(\phi,\theta,\varepsilon)\cos \phi \cos^{2} \theta \nonumber \\ +
\varepsilon b(\phi,\theta,\varepsilon)\sin \theta \cos \theta \cos \phi.
\end{eqnarray}
If $ \sin \theta \neq 0$, we divide the first equation in (\ref{E:EqDiffNew15}) by $\sin \phi \sin \theta $, to get
\begin{equation}\label{E:EqDiffNew161}
\dot{\theta}= -\left [ \left | X_{0} \right | + \varepsilon a(\phi, \theta, \varepsilon)\right ]+
\varepsilon b(\phi, \theta, \varepsilon)\sin \theta \cot \phi+ \\
\varepsilon c(\phi, \theta, \varepsilon) \cos \theta \cot \phi.
\end{equation}
Thus, the second differential equation of (\ref{E:EqDiffNew15}) is satisfied.\\
If $ \sin \theta = 0$, then $\cos \theta \neq 0$ and we divide the second equation in (\ref{E:EqDiffNew15}) by
$\sin \phi \cos \theta$ to get
\begin{equation}\label{E:EqDiffNew16}
\dot{\theta}= -\left [\left | X_{0} \right | + \varepsilon a(\phi, \theta, \varepsilon)\right ]+
\varepsilon b(\phi, \theta, \varepsilon)\sin \theta \cot \phi+\\
\varepsilon c(\phi, \theta, \varepsilon) \cos \theta \cot \phi.
\end{equation}
Thus, the first of the differential equations (\ref{E:EqDiffNew15}) is satisfied.\\
Therefore, we obtain the system
\begin{equation}\label{E:EqDiffNew17}
\begin{array}{lll}
\dot{\phi} &=& -\varepsilon b(\phi, \theta, \varepsilon) \cos \theta +
\varepsilon c(\phi, \theta, \varepsilon) \sin \theta, \\
\dot{\theta} &=& -\left [\left | X_{0} \right | + \varepsilon a(\phi,\theta, \varepsilon)\right ]+
\varepsilon b(\phi, \theta, \varepsilon)\sin \theta \cot \phi+\varepsilon c(\phi, \theta, \varepsilon)
\cos \theta \cot \phi.\\
\end{array}
\end{equation}
or
\begin{equation}\label{E:EqDiffNew5}
\begin{array}{lrl}
\dot \phi &=& \varepsilon \left [-b(\phi, \theta, \varepsilon)\cos \theta+c(\phi, \theta, \varepsilon)
\sin \theta \right ],\\
\dot \theta &=& -\left | X_{0} \right | + \varepsilon \{a(\phi, \theta, \varepsilon)+
\left [b(\phi, \theta, \varepsilon)\sin \theta + c(\phi, \theta, \varepsilon)\cos \theta \right ]\cot \phi\}.\\
\end{array}
\end{equation}
For $\phi \in (0, \pi)$, $\theta \in [0,2\pi)$ and for $\varepsilon >0$ small enough, we can have that
\[ -\left | X_{0} \right | + \varepsilon \{a(\phi, \theta, \varepsilon)+\left [b(\phi, \theta, \varepsilon)\sin \theta +
c(\phi, \theta, \varepsilon)\cos \theta \right ]\cot \phi\} > 0\] and we can rescale time, choosing $\tau= \tau(t)$
such that
\begin{equation}\label{E:tau}
\frac{d \tau}{d t}= -\left | X_{0} \right | + \varepsilon \left \{a(\phi, \theta, \varepsilon)+
\left [b(\phi, \theta, \varepsilon)\sin \theta + c(\phi, \theta, \varepsilon)\cos \theta \right ]\cot \phi \right \}.
\end{equation}
We get
\begin{equation}\label{E:rescaled_phi}
\begin{array}{lll}
\frac{d\phi}{d\tau} &=& \frac{d\phi}{d t} \frac{d t}{d \tau}=
\frac{\varepsilon \left [-b(\phi, \theta, \varepsilon)\cos \theta+
c(\phi, \theta, \varepsilon)\sin \theta \right ]}{-\left | X_{0} \right | +
\varepsilon \left \{a(\phi, \theta, \varepsilon)+
\left [b(\phi, \theta, \varepsilon)\sin \theta + c(\phi, \theta, \varepsilon)\cos \theta \right ]\cot \phi \right \}},\\
\frac{d \theta}{d \tau} &=& \frac{d \theta}{d t}\frac{d t}{d \tau}=1.\\
\end{array}
\end{equation}
Since \[ b(\phi, \theta, \varepsilon)=  b(\phi, \theta, 0)+ \varepsilon h_{1}(\phi, \theta, \varepsilon)\] and
\[ c( \phi, \theta, \varepsilon)=  c( \phi, \theta, 0)+ \varepsilon h_{2}( \phi, \theta, \varepsilon), \]
it results that
\begin{equation}\label{E:EqDiffNew6}
\begin{array}{lll}
\frac{d\phi}{d\tau} &=& \frac{\varepsilon}{\left | X _{0} \right | }\left [ b(\phi, \theta, 0)\cos \theta-
c(\phi, \theta, 0)\sin \theta \right ]+\varepsilon^{2} h(\phi, \theta, \varepsilon),\\
\frac{d\theta}{d\tau} &=& 1.\\
\end{array}
\end{equation}
Let us construct the Poincar\'{e} map ( see \cite{Wi}, that is the time $2\pi$ map) for the flow given by the
differential equations (\ref{E:EqDiffNew6}). The Poincar\'{e} section associated to the flow given by the
differential equations (\ref{E:EqDiffNew6}) is given by $\theta= 0$, $\phi \in (0, \pi)$, with $\phi$ chosen
independent of $\theta_{0} \in (0,\pi) $. The Poincar\'{e} map is given by
\begin{equation*}
P(\phi_{1}, \varepsilon)= \phi(2\pi, \phi_{1}, \varepsilon),
\end{equation*}
where $\phi(t, \phi_{1}, \varepsilon)$ is the solution of the following initial problem after we relabel $\tau=t$,
\begin{equation}\label{E:PoincareIVP}
\begin{array}{lll}
\frac{d\phi}{d t} &=& \frac{\varepsilon}{\left | X _{0} \right | }\left [b(\phi, t, 0)\cos t-
c(\phi, t, 0)\sin t \right ]+ \varepsilon^{2} h(\phi, t, \varepsilon),\\
\phi(0) &=& \phi_{1}.\\
\end{array}
\end{equation}
Therefore, we get by using (\ref{E:PoincareIVP})
\begin{equation*}
\begin{split}
P(\phi_{1}, \varepsilon) &= \phi(0,\phi_{1},\varepsilon)+
\int_{0}^{2\pi}\frac{d \phi}{d t}(t, \phi_{1}, \varepsilon) \, dt\\
&=\phi_{1}+ \int_{0}^{2\pi} \left \{\frac{\varepsilon}{\left | X_{0}
\right | }\left [b(\phi(t,\phi_{1},\varepsilon), t, 0)\cos t-
c(\phi(t, \phi_{1}, \varepsilon), t, 0)\sin t \right ]\right .\\
& \left . +\varepsilon^{2} h(\phi(t,\phi_{1},\varepsilon), t, \varepsilon)
\right \} \, dt
\end{split}
\end{equation*}
or
\begin{equation*}
P(\phi_{1}, \varepsilon)=\phi_{1}+ \int_{0}^{2\pi} \frac{\varepsilon}{\left | X_{0} \right | }
\left [ b(t,\phi_{1},0)\cos t- c(t,\phi_{1}, 0)\sin t \right ] \, dt +\varepsilon^{2}h_{*}(\phi_{1}, \varepsilon).
\end{equation*}
We have that
\begin{equation*}
P(\phi_{1}, \varepsilon)= \phi_{1}+\frac{\varepsilon}{\left | X_{0} \right | }I(\phi_{1})+
\varepsilon^{2} h_{*}(\phi_{1}, \varepsilon),
\end{equation*}
where we define
\begin{equation*}
I(\phi_{1})= \int_{0}^{2\pi} \left [b(\phi_{1},t,0)\cos t-c(\phi_{1},t,0)\sin t \right ] \, dt.
\end{equation*}
Then, \[P(\phi_{0},0)-\phi_{0}= 0,\] since $x(t,x_{0})= e^{-X_{0}t}x_{0}$ is a periodic solution of the differential
equations (\ref{E:pertdif_projected}) for $\varepsilon=0$, that is $s(t,s_{0})= e^{-X_{00}t}s_{0}$ is a periodic
solution of the differential equations (\ref{E:EqDiffNew2}) for $\varepsilon=0$ and
$s_{0}= \left ( \begin{array}{c}
      \sin \phi_{0}\cos \theta_{0} \\
      \sin \phi_{0} \sin \theta_{0}\\
      \cos \phi_{0}\\
                \end{array} \right )$.
Let us consider the equation
\begin{equation*}
P(\phi_{1},\varepsilon)-\phi_{1}= \frac{\varepsilon}{\left | X_{0} \right |}\left [I(\phi_{1})+
\varepsilon h_{*}(\phi_{1}, \varepsilon) \right ]= 0,
\end{equation*}
that is
\begin{equation*}
I(\phi_{1})+\varepsilon h_{*}(\phi_{1}, \varepsilon)= 0.
\end{equation*}
If
\begin{equation*}
I(\phi_{0})= 0 \mbox { and } I^{'}(\phi_{0}) \neq 0,
\end{equation*}
then, using the implicit function theorem, we find for $\varepsilon \geq 0$ small a sufficiently smooth
branch $\phi(\varepsilon)$ of fixed points of $P$ such that $\phi(0)= \phi_{0}$.\\
We have \[b(\phi_{1},t,0)= G_{2}(\sin \phi_{1} \cos t, \sin \phi_{1} \sin t, \cos \phi_{1},0)\] and
\[c(\phi_{1},t,0)= G_{3}(\sin \phi_{1} \cos t, \sin \phi_{1} \sin t, \cos \phi_{1},0).\]\\
Therefore,
\begin{equation*}
\begin{split}
\frac{\partial{b}}{\partial{\phi_{1}}}(\phi_{1},t,0) &=
\frac{\partial{G_{2}}}{\partial{x_{1}}}(\sin \phi_{1} \cos t, \sin \phi_{1} \sin t, \cos \phi_{1},0)\cos \phi_{1}\cos t \\ &+
\frac{\partial{G_{2}}}{\partial{x_{2}}}(\sin \phi_{1} \cos t, \sin \phi_{1} \sin t, \cos \phi_{1},0)\cos \phi_{1} \sin t \\&-
\frac{\partial{G_{2}}}{\partial{x_{3}}}(\sin \phi_{1} \cos t, \sin \phi_{1} \sin t, \cos \phi_{1},0)\sin \phi_{1}
\end{split}
\end{equation*}
and
\begin{equation*}
\begin{split}
\frac{\partial{c}}{\partial{\phi_{1}}}(\phi_{1},t,0) &=
\frac{\partial{G_{3}}}{\partial{x_{1}}}(\sin \phi_{1} \cos t, \sin \phi_{1} \sin t, \cos \phi_{1},0)\cos \phi_{1}\cos t \\ &+
\frac{\partial{G_{3}}}{\partial{x_{2}}}(\sin \phi_{1} \cos t, \sin \phi_{1} \sin t, \cos \phi_{1},0)\cos \phi_{1} \sin t \\ &-
\frac{\partial{G_{3}}}{\partial{x_{3}}}(\sin \phi_{1} \cos t, \sin \phi_{1} \sin t, \cos \phi_{1},0)\sin \phi_{1}.
\end{split}
\end{equation*}
For $\varepsilon \geq 0$ small, we get a sufficiently smooth branch $s^{\varepsilon}(t)$ of periodic solutions for the
differential equations (\ref{E:EqDiffNew2}) such that $s^{0}(t)= e^{-X_{00}t}s_{0}$.\\
Thus, for $\varepsilon \geq 0$ small, we get a sufficiently smooth branch $x^{\varepsilon}(t)$ of periodic solutions for
the differential equations (\ref{E:pertdif_projected}) such that $x^{0}(t)= e^{-X_{0}t}x_{0}$.\\
If $I^{'}(\phi_{0}) < 0 $, the fixed point $\phi(\varepsilon)$ is locally asymptotically stable for $P$, and if
$I^{'}(\phi_{0}) > 0$, the fixed point $\phi(\varepsilon)$ is unstable for $P$.\\
The stability of the periodic solutions $x^{\varepsilon}(t)$ which persist for $\varepsilon >0$ is the same as the
stability of the fixed points $\phi(\varepsilon)$ of the Poincar\'{e} map $P$ from which the periodic solutions
$x^{\varepsilon}(t)$ are obtained.\\
This ends the proof of the second conclusion of Theorem \ref{thm:FSB_DE_Orbit_E}.
\end{proof}
\begin{proof}[Proof of Proposition \ref{prop:FSB_Phase_E}]
We prove the first conclusion for the equilibria $x^{1}(\varepsilon)$. A similar proof can be done for the equilibria
$x^{2}(\varepsilon)$. \\
Let $A_{\varepsilon} \in SO(3)$ be such that $x_{1}^{Q}= A_{\varepsilon}x^{1}(\varepsilon)$. We can choose a
sufficiently smooth branch $A_{\varepsilon}$, by taking $C_{\varepsilon}=(\pi_{s} \circ \beta^{-1})(x^{1}(\varepsilon))$
and then $A_{\varepsilon}=C_{\varepsilon}^{-1}$. We have $C_{\varepsilon}x_{1}^{Q}=x^{1}(\varepsilon)$.\\
Taking into account that $\mathbf{\Psi}(t,C,\varepsilon)=[\mathbf{\Phi}(t,A,\varepsilon)]^{-1}$,
with $A=C^{-1}$, we check that the differential equations (\ref{E:pertdif_changed}) have a solution of the form
\begin{equation*}
C^{RW}(t, \varepsilon) \overset{\text{def}}= C_{\varepsilon}e^{-\alpha_{1}(\varepsilon)Qt} \mbox{ with }
\alpha_{1}(\varepsilon)= \left | X_{0} \right | +O(\varepsilon).
\end{equation*}
By the definitions of the projected flows $\widetilde{\mathbf{\Psi}_{1}}$ and $\widetilde{\mathbf{\Psi}}$,
we have $\widetilde{\mathbf{\Psi}_{1}}(t,C_{\varepsilon}x_{1}^{Q},\varepsilon) =C_{\varepsilon}x_{1}^{Q}=
x^{1}(\varepsilon)$, then
\begin{equation*}
\widetilde{\mathbf{\Psi}}(t, C_{\varepsilon}\cdot SO(2),\varepsilon)=
\pi(\mathbf{\Psi}(t, C_{\varepsilon},\varepsilon))=  C_{\varepsilon} \cdot SO(2),
\end{equation*}
which implies
\begin{equation}\label{E:PhaseRW1}
C^{RW}(t, \varepsilon)= \mathbf{\Psi}(t, C_{\varepsilon},\varepsilon)= C_{\varepsilon}e^{\tau(t,\varepsilon)Q},
\end{equation}
where $\tau(0, \varepsilon)= 0$. Without loss of generality, we assume that $\tau(t, \varepsilon)$ sufficiently smooth
with respect to $t$. It follows that $C^{RW}(t, \varepsilon)$ verifies the differential equations (\ref{E:pertdif_changed}).\\
If we substitute (\ref{E:PhaseRW1}) into the differential equations (\ref{E:pertdif_changed}), we get
\begin{equation*}
C_{\varepsilon}\frac{\partial{\tau}}{\partial{t}}(t, \varepsilon)Q e^{\tau(t,\varepsilon)Q}=
-\left [X_{0}+\varepsilon \widetilde{g}(C_{\varepsilon}e^{\tau(t,\varepsilon)Q}, \varepsilon)\right ]
C_{\varepsilon}e^{\tau(t,\varepsilon)Q}
\end{equation*}
or using the $SO(2)$-invariance of $\widetilde{g}(.,\varepsilon)$, we get
\begin{equation*}
\begin{split}
C_{\varepsilon}\frac{\partial{\tau}}{\partial{t}}(t, \varepsilon)Q &=
-\left [X_{0}+ \varepsilon \widetilde{g}(C_{\varepsilon}, \varepsilon)\right ]C_{\varepsilon} \\ &
\Rightarrow \frac{\partial{\tau}}{\partial{t}}(t, \varepsilon)Q=
-C_{\varepsilon}^{-1}\left [X_{0}+ \varepsilon \widetilde{g}(C_{\varepsilon}, \varepsilon)\right ]C_{\varepsilon}\\ &
\Rightarrow \left | \frac{\partial{\tau}}{\partial{t}}(t,\varepsilon)Q \right |= \left | -
C_{\varepsilon}^{-1}\left [X_{0}+\varepsilon \widetilde{g}(C_{\varepsilon}, \varepsilon)\right ]C_{\varepsilon} \right |
\end{split}
\end{equation*}
or, using $\left | Q \right | =1$, we get
\begin{equation}
\begin{split}
\left | \frac{\partial{\tau}}{\partial{t}}(t, \varepsilon) \right | &=
\left | X_{0}+ \varepsilon \widetilde{g}(C_{\varepsilon}, \varepsilon)\right | \\ &
\Rightarrow \frac{\partial{\tau}}{\partial{t}}(t, \varepsilon)= \mbox{ constant }
\overset{\text{def}}= \alpha_{1}(\varepsilon).
\end{split}
\end{equation}
Therefore,
\begin{equation*}
\tau(t, \varepsilon)= -\alpha_{1}(\varepsilon)t+\tau(0,\varepsilon)\Rightarrow \tau(t, \varepsilon)=
-\alpha_{1}(\varepsilon)t.
\end{equation*}
Also, we get $\alpha_{1}(\varepsilon)Q= C_{\varepsilon}^{-1}\left [X_{0}+
\varepsilon \widetilde{g}(C_{\varepsilon},\varepsilon)C_{\varepsilon}\right ]$ and $\alpha_{1}(\varepsilon)$ is
sufficiently smooth. Then, \[\alpha_{1}(\varepsilon)=\left | X_{0} \right | +O(\varepsilon),\]
since $\alpha_{1}(0)Q=C_{0}^{-1}X_{0}C_{0}$ implies $\alpha_{1}(0)\overrightarrow Q=
\left | X_{0} \right | \overrightarrow Q $, namely $\alpha_{1}(0)=\left | X_{0} \right |$
because $x_{1}^{Q}=C_{0}^{-1}x^{1}(0)$, $x^{1}(0)=x^{0}_{1}$.\\
The stability issue for $C^{RW}(t,\varepsilon)$ is proved in \cite{Ch}. It is now easy to get the results for the
differential equations (\ref{E:pertdif}) using the fact that $\mathbf{\Psi}(t,C,\varepsilon)=
[\mathbf{\Phi}(t,A,\varepsilon)]^{-1}$, where $A=C^{-1}$.\\
This proves the first conclusion. Similarly, we can prove second conclusion of Proposition \ref{prop:FSB_Phase_E}.\\
We now prove the conclusion (3) of Proposition \ref{prop:FSB_Phase_E}.\\
Let $D(t, \varepsilon) \in SO(3)$ be such that $D(t, \varepsilon)x_{1}^{Q}= x^{\varepsilon}(t)$.
Let $T(\varepsilon)= \frac{2\pi}{\left | X_{0} \right |}+O(\varepsilon)$ be the period of the function
$x^{\varepsilon}(t)$. We can choose a sufficiently smooth  branch $D(t, \varepsilon)$, by taking
$D(t,\varepsilon)= (\pi_{s} \circ \beta^{-1})(x^{\varepsilon}(t))$. Since $x^{\varepsilon}(t)$ is
$T(\varepsilon)$-periodic, we have $D(t+T(\varepsilon),\varepsilon)x_{1}^{Q}=D(t,\varepsilon)x_{1}^{Q}$.
Therefore, since the section $\pi_{s}$ is local, we can only say that there exists a sufficiently smooth function
$g(t,\varepsilon)$ such that
\begin{equation}\label{E:dd}
D(t+T(\varepsilon),\varepsilon)=D(t,\varepsilon)e^{g(t,\varepsilon)Q}.
\end{equation}
Without loss of generality, we assume that $g(0,\varepsilon)=0$, otherwise we take
$D_{1}(t,\varepsilon)= D(t,\varepsilon)e^{-g(0,\varepsilon)Q}$.\\\\
Taking into the account that $\mathbf{\Psi}(t,C,\varepsilon)=[\mathbf{\Phi}(t,A,\varepsilon)]^{-1}$,
where $A=C^{-1}$, we check that the differential equations (\ref{E:pertdif_changed}) have a solution of the
form \[C^{MRW}(t, \varepsilon)= B^{*}(t, \varepsilon)e^{-\beta(\varepsilon)Qt}\] with
$\beta(\varepsilon)= O(\varepsilon)$, $B^{*}(0, \varepsilon)= D(0, \varepsilon)$ and $B^{*}(t, \varepsilon)$ is
$T(\varepsilon)$-periodic.\\\\
Let $C_{\varepsilon}^{1}= D(0, \varepsilon)$. We construct $B^{*}(t,\varepsilon)$ and $\beta(\varepsilon)$.\\
By the definitions of the projected flow $\widetilde{\mathbf{\Psi}_{1}}$ and $\widetilde{\mathbf{\Psi}}$,
we have $\widetilde{\mathbf{\Psi}_{1}}(t,C_{\varepsilon}^{1}x_{1}^{Q},\varepsilon)=
x^{\varepsilon}(t)=D(t,\varepsilon)x_{1}^{Q}$ or $\mathbf{\Psi}(t,C_{\varepsilon}^{1} \cdot
SO(2),\varepsilon)=\pi(\mathbf{\Psi}(t,C_{\varepsilon}^{1},\varepsilon))=D(t,\varepsilon) \cdot SO(2)$, which implies
\begin{equation}\label{E:PhaseMRW2}
\mathbf{\Psi}(t, C^{1}_{\varepsilon}, \varepsilon)= D(t,\varepsilon)e^{-\tau(t, \varepsilon)Q}
\mbox{ with }\tau(0, \varepsilon)=0.
\end{equation}
Without loss of generality, we assume that $\tau(t, \varepsilon)$ is a sufficiently smooth function with respect to $t$.
If we substitute $(\ref{E:PhaseMRW2})$ into the differential equations (\ref{E:pertdif_changed}), we get
\begin{equation*}
\dot{D}=-\left [X_{0}+\varepsilon \widetilde{g}(D,\varepsilon)\right ]D+\dot{\tau}D Q \Rightarrow \dot{\tau}D Q=
D^{-1}\dot{D}+D^{-1} \left [X_{0}+\varepsilon \widetilde{g}(D,\varepsilon)\right ]D.
\end{equation*}
Therefore, the function $\tau(t, \varepsilon)$ is sufficiently smooth in $\varepsilon$ and in $t$. We define
\begin{equation*}
C^{MRW}(t,\varepsilon) \overset{\text{def}} = D(t,\varepsilon)e^{-\tau(t, \varepsilon)Q}
e^{\beta(\varepsilon)Qt}e^{-\beta(\varepsilon)Qt}, \mbox { where } \beta(\varepsilon)=
\frac{\tau(T(\varepsilon),\varepsilon)}{T(\varepsilon)}.
\end{equation*}
Let us define $B^{*}(t,\varepsilon)= D(t,\varepsilon)e^{-\tau(t, \varepsilon)Q}e^{\beta(\varepsilon)Qt}$.
Then, $C^{MRW}(t,\varepsilon)= B^{*}(t,\varepsilon)e^{-\beta(\varepsilon)Qt}$.\\\\
We have that $B^{*}(T(\varepsilon),\varepsilon)= D(T(\varepsilon),\varepsilon)e^{-\tau(T(\varepsilon),\varepsilon)Q}
e^{\beta(\varepsilon)T(\varepsilon)Q}=D(0,\varepsilon)e^{g(0,\varepsilon)Q}= D(0,\varepsilon)=B^{*}(0,\varepsilon)$
by the definition of $\beta(\varepsilon)$, relation (\ref{E:dd}) and the fact that $g(0,\varepsilon)=0$.\\
Since $\tau$ is sufficiently smooth in $t$ and $\varepsilon$ and $T(\varepsilon)$ is sufficiently smooth,
it follows from the definition of $\beta(\varepsilon)$ that $\beta(\varepsilon)$ is a sufficiently smooth function and
$\beta(\varepsilon)=  \frac{\tau(T(0),0)}{T(0)}+O(\varepsilon)$.\\\\
We check that $B^{*}(t,\varepsilon)$ is $T(\varepsilon)$-periodic. If we substitute $C= Be^{-\beta(\varepsilon)Qt}$ into
the differential equations (\ref{E:pertdif_changed}), we get
\begin{equation*}
\dot{B}e^{-\beta(\varepsilon)Qt}+B(-\beta(\varepsilon)Q)e^{-\beta(\varepsilon)Qt}=
-\left [ X_{0}+\varepsilon \widetilde{g}(B,\varepsilon)\right ] B e^{-\beta(\varepsilon)Qt}
\end{equation*}
or
\begin{equation}\label{E:eq112}
\dot{B}= \beta (\varepsilon) B Q- \left [ X_{0}+\varepsilon \widetilde{g}(B,\varepsilon) \right ] B .
\end{equation}
Since $C^{MRW}(t,\varepsilon)= B^{*}(t,\varepsilon)e^{-\beta(\varepsilon)Qt}$ and $C^{MRW}(t,\varepsilon)$
is a solution of the differential equations (\ref{E:pertdif_changed}), we get that $ B^{*}(t,\varepsilon)$ is a
solution of the differential equations (\ref{E:eq112}). Since $B^{*}(T(\varepsilon),\varepsilon)= B^{*}(0,\varepsilon)$,
we get that $C^{*}(t,\varepsilon)= B^{*}(t+T(\varepsilon),\varepsilon)$ is also a solution of the differential equations
(\ref{E:eq112}) such that $C^{*}(0,\varepsilon)=B^{*}(0,\varepsilon)$. Therefore, the function
$B^{*}(t,\varepsilon)$ is $T(\varepsilon)$-periodic.\\\\
We check that $\beta(\varepsilon)=O(\varepsilon)$ up to $k\frac{2\pi}{T(\varepsilon)}$, for some $k \in \mathbb{Z}$.
Therefore, we show that $\tau(T(0),0)=2k\pi$ for some $k \in \mathbb{Z}$.\\
By the definition of $D$ we get that $D(t,0)x_{1}^{Q}=x^{0}(t)= x(t,x_{0})e^{-X_{0}t}x_{0}$ or
\begin{equation}\label{E:abcd}
D(t,0)=e^{-X_{0}t}D(0,0)e^{g_{1}(t)Q},
\end{equation}
with $g_{1}$ sufficiently smooth such that $g_{1}(0)=0$. \\\\
We show that $g_{1}(T(0))=2k\pi$ for some $k \in \mathbb{Z}$ and $\tau(T(0),0)=g_{1}(T(0))$.\\
Using (\ref{E:dd}), it follows that $D(t+T(0),0)= D(t,0)e^{g(t,0)Q}$. Since $g(0,0)=0$, we get
$D(T(0),0)=D(0,0)=e^{-X_{0}T(0)}D(0,0)e^{g_{1}(T(0))Q}$. Taking into account that $\left | X(0) \right | T(0) =2\pi$,
we get $e^{g_{1}(T(0))Q}=I_{3}$. It follows that $g_{1}(T(0))=2k\pi$ for some $k \in \mathbb{Z}$.\\
Also, $\mathbf{\Psi}(t, C^{1}_{0}, 0)= D(t, 0)e^{-\tau(t,0)Q}= e^{-X_{0}t}D(0,0)e^{g_{1}(t)Q}e^{-\tau(t,0)Q}$ by
(\ref{E:abcd}). On the other hand, $\mathbf{\Psi}(t, C^{1}_{0}, 0)=e^{-X_{0}t}D(0,0)$. Then, we
get $\tau(t,0)=g_{1}(t)+2l(t)\pi$ with $l(t) \in \mathbb{Z}$. Since $g_{1}$ and $\tau$ are sufficiently smooth,
we get that $l(t)$ is sufficiently smooth and, since $l(t) \in \mathbb{Z}$, we get $l(t)=$constant. Since
$\tau(0,0)=g_{1}(0)=0$, we get $l=0$.\\
Thus $\tau(t,0)=g_{1}(t)$ implies $\tau(T(0),0)=g_{1}(T(0))=2k\pi$.\\\\
Then, $\beta(0)= \frac{\tau(T(0),0)}{T(0)}=\frac{g_{1}(T(0))}{T(0)}= \frac{2k\pi}{T(0)}$.\\\\
Therefore, $\beta(\varepsilon)= \frac{2k\pi}{T(0)}+O(\varepsilon)= k\frac{2\pi}{T(\varepsilon)}+\beta_{1}(\varepsilon)$,
where $\beta_{1}(\varepsilon)=O(\varepsilon)$.\\
We have $C^{MRW}(t,\varepsilon)= B^{*}(t,\varepsilon)e^{-\beta(\varepsilon)Qt}=
B^{*}(t,\varepsilon)e^{-\frac{2\pi}{T(\varepsilon)}kQt}e^{-\beta_{1}(\varepsilon)Qt}$.\\
Let $B_{1}^{*}(t,\varepsilon)=B^{*}(t,\varepsilon)e^{-\frac{2\pi}{T(\varepsilon)}kQt}$,
$B_{1}^{*}(t, \varepsilon)$ is $T(\varepsilon)$-periodic and $B_{1}^{*}(0,\varepsilon)= B^{*}(0,\varepsilon)$.\\
We may drop the subscript 1 in $B_{1}^{*}$ and $\beta_{1}(\varepsilon)$.\\
The stability issue for $C^{MRW}(t,\varepsilon)$ is proved in \cite{Ch}. It is now easy to get the results for the
differential equations (\ref{E:pertdif}) using the fact that $\mathbf{\Psi}(C,t)=[\mathbf{\Phi}(A,t)]^{-1}$, where $A=C^{-1}$.
This ends the proof of the third conclusion of Proposition \ref{prop:FSB_Phase_E} and its proof.
\end{proof}
\begin{proof}[Proof of Theorem \ref{thm:results_FSB}]
It is a consequence of Theorem \ref{thm:CMR_FSB_RE}, Proposition \ref{prop:FSB_Phase_E} and of the $SO(2)$-equivariance
of the flow $\mathbf{\Phi}$.\\
\end{proof}

\end{document}